\RequirePackage{fix-cm}
\documentclass[smallextended]{svjour3}       
\smartqed  
\usepackage{graphicx}
\usepackage{amsmath}
\usepackage{amsfonts}
\usepackage{amssymb}
\usepackage{array,longtable}
\usepackage{amscd}
\usepackage[cp1251]{inputenc}
\usepackage[english]{babel}
\usepackage[usenames]{color}
\usepackage{graphicx}
\graphicspath{{pictures/}}
\DeclareGraphicsExtensions{.pdf,.png,.jpg}
\usepackage[matrix,arrow,curve]{xy}
\begin{document}
\setcounter{tocdepth}{5}
\title{Generalized $D$-stability and diagonal dominance with applications to stability and transient response properties of systems of ODE 
}

\titlerunning{Generalized $D$-stability and diagonal dominance}        

\author{Olga Y. Kushel         \and
        Raffaella Pavani
}


\institute{O. Kushel \at
              Shanghai University, \\ Department of Mathematics, \\ Shangda Road 99, \\ 200444 Shanghai, China \\
              Tel.: +86(021)66133292\\
              \email{kushel@mail.ru}           
           \and
           R. Pavani \at
              Politecnico di Milano, \\ Department of Mathematics, \\ Piazza L. da Vinci 32, \\ 20133 Milano, Italia \\
              Tel.: +39(02)23994516\\
              \email{raffaella.pavani@polimi.it}
}

\date{Received: date / Accepted: date}

\maketitle

\begin{abstract}
In this paper, we introduce the class of diagonally dominant (with respect to a given LMI region ${\mathfrak D} \subset {\mathbb C}$) matrices that possesses the analogues of well-known properties of (classical) diagonally dominant matrices, e.g their spectra are localized inside the region $\mathfrak D$. Moreover, we show that in some cases, diagonal $\mathfrak D$-dominance implies $({\mathfrak D}, {\mathcal D})$-stability ( i.e. the preservation of matrix spectra localization under multiplication by a positive diagonal matrix). Basing on the properties of diagonal stability and diagonal dominance, we analyze the conditions for stability of second-order dynamical systems. We show that these conditions are preserved under system perturbations of a specific form (so-called $D$-stability). We apply the concept of diagonal $\mathfrak D$-dominance to the analysis of the minimal decay rate of second-order systems and its persistence under specific perturbations (so-called relative $D$-stability). Diagonal $\mathfrak D$-dominance with respect to some conic region $\mathfrak D$ is also shown to be a sufficient condition for stability and $D$-stability of fractional-order systems.

\end{abstract}

\section{Introduction}
\subsection{Notations}
Here, we introduce the following notations:
\begin{enumerate}
\item[] ${\mathcal M}^{n \times n}$ for the set of real $n \times n$ matrices;
\item[] ${\mathcal M}^{n \times n}({\mathbb C})$ for the set of complex $n \times n$ matrices;
\item[] $\sigma (\mathbf{A})$ for the
spectrum of a matrix $\mathbf{A}$ $\in \mathcal{M}^{n\times n}$ (i.e. the set of all eigenvalues of $\mathbf A$ defined as zeroes of its characteristic polynomial $f_{\mathbf A}(\lambda):= \det(\lambda{\mathbf I - {\mathbf A}})$);
\item[] ${\mathbb C}^-$ for the open left-hand side of the complex plane $\mathbb C$, i.e. $${\mathbb C}^- = \{z \in {\mathbb C}: {\rm Re}(z) < 0\};$$
\item[] $\overline{{\mathbb C}^+}$ for the closed right-hand side of the complex plane $\mathbb C$, i.e. $$\overline{{\mathbb C}^+} = \{z \in {\mathbb C}: {\rm Re}(z) \geq 0\};$$
\item[] ${\mathcal D}^+ \subset \mathcal{M}^{n\times n}$ for the set of all positive diagonal matrices (i.e. matrices with
positive entries on the principal diagonal while the entries outside the principal diagonal are zeroes);
\item[] ${\mathcal D}^+_{(0,1]}\subset \mathcal{M}^{n\times n}$ for the subclass of ${\mathcal D}^+$, defined as follows:
$${\mathcal D}^+_{(0,1]} = \{{\mathbf D} = {\rm diag}\{d_{11}, \ \ldots, \ d_{nn}\}: 0 < d_{ii} \leq 1, \ i = 1, \ \ldots, \ n\}; $$
\item[] ${\mathcal D}^+_{\geq 1}\subset \mathcal{M}^{n\times n}$ for the subclass of ${\mathcal D}^+$, defined as follows:
$${\mathcal D}^+_{\geq 1} = \{{\mathbf D} = {\rm diag}\{d_{11}, \ \ldots, \ d_{nn}\}: 1 \leq d_{ii} < + \infty, \ i = 1, \ \ldots, \ n\}. $$
\end{enumerate}

\subsection{Diagonal dominance and matrix stability} In this paper, we generalize the concept of diagonal dominance. First, let us recall the classical definition of diagonally dominant matrices (see, for example, \cite{HOJ}, also \cite{NIS}, \cite{NIS2}).

{\bf Definition 1.} A matrix ${\mathbf A} \in {\mathcal M}^{n \times n}({\mathbb C})$ is called {\it strictly row diagonally dominant} if the following inequalities hold:
 \begin{equation}\label{ROW}|a_{ii}| > \sum_{i\neq j}|a_{ij}| \qquad i = 1, \ \ldots, \ n.\end{equation}
 A matrix $\mathbf A$ is called {\it strictly column diagonally dominant} if ${\mathbf A}^T$ is strictly row diagonally dominant.

{\bf Definition 1'.} A matrix ${\mathbf A} \in {\mathcal M}^{n \times n}({\mathbb C})$ is called {\it generalized diagonally dominant} if there exist positive scalars (weights) $m_i$, $i = 1, \ \ldots, \ n$ such that
$$m_i|a_{ii}|>\sum_{i\neq j }m_j|a_{ij}|, \qquad i = 1, \ \ldots, \ n. $$
(i.e. if there is a positive diagonal matrix ${\mathbf M} = {\rm diag}\{m_1, \ \ldots, \ m_n\}$, such that ${\mathbf A}{\mathbf M}$ is strictly row diagonally dominant).
If, in addition, $a_{ii} < 0$, $i = 1, \ \ldots, \ n,$ then ${\mathbf A}$ is called {\it negative diagonally dominant (NDD)}.

It is well-known (see, for example, \cite{TAK}, p. 382, Theorem 4.C.2) that if a complex matrix ${\mathbf A}$ is NDD, then $\mathbf A$ is Hurwitz stable, i.e. every $\lambda \in \sigma (\mathbf{A})$ satisfies ${\rm Re}(\lambda) < 0$. The stability of $\mathbf A$ implies the Lyapunov asymptotic stability of a system of ODE with the system matrix $\mathbf A$. Being a sufficient for stability condition, negative diagonal dominance is of particular interest by itself due to its connection to the properties of nonlinear systems (see \cite{SIL}). Less known is the fact that negative diagonal dominance implies some special concepts of matrix and system stability, that are stronger than just stability and shows stability preservation under specific matrix (system) perturbations. Such concepts include multiplicative $D$-stability (see \cite{JOHN1}, \cite{KEK}) and diagonal stability (see \cite{MOY}).

{\bf Definition 2}. A matrix $\mathbf{A}$ $\in \mathcal{M}^{n\times n}$ is
called {\it (multiplicative) $D$-stable} if ${\rm Re}(\lambda )<0$ for all $\lambda \in
\sigma (\mathbf{DA})$, where ${\mathbf D}$ is any matrix from ${\mathcal D}^+$.

{\bf Definition 3}. A matrix ${\mathbf A} \in {\mathcal M}^{n \times n}$ is called {\it diagonally stable} if the matrix $${\mathbf W}:={\mathbf D}{\mathbf A} + {\mathbf A}^{T}{\mathbf D}$$
is negative definite for some positive diagonal matrix ${\mathbf D}$.

The above matrix classes play an important role in the theory of stability (see, for example, \cite{LOG}, \cite{KAB}, \cite{KU2}) and have numerous applications in the economics, mathematical ecology, mechanics and other branches of science. The proof of the fact that NDD matrices are $D$-stable is based on the simple observance that the negative diagonal dominance implies stability and is preserved under multiplication by a positive diagonal matrix (see \cite{JOHN1}, p. 54, Observation (i)). Note, that the property of diagonal dominance can be proved (disproved) by a finite number of steps while the conditions of $D$-stability of $\mathbf A$ involves checking all the products of the form ${\mathbf D}{\mathbf A}$, where ${\mathbf D}$ runs along the infinite set of positive diagonal matrices.

The concept of negative diagonal dominance is also applied to establishing the stability of systems of second-order differential equations (see \cite{NIS}, \cite{NIS2}), with the applications to the stability of mechanical systems and certain economic models.
\subsection{Paper overview} In this paper, we are concerned with the problem of a matrix spectra localization inside a prescribed convex region ${\mathfrak D} \subset {\mathbb C}$ (so called $\mathfrak D$-stability problem) and its robust aspects. Many problems of the system dynamics lead to establishing ${\mathfrak D}$-stability of the system matrix with respect to a specified region ${\mathfrak D} \subset {\mathbb C}$ (see, for example, \cite{GUT2}, \cite{GUJU}). The regions of particular interest are LMI (Linear Matrix Inequality) regions, that include the shifted left half-plane, the unit disk and the conic sector around the negative direction of the real axis. Due to the rapid advance in control theory and its application, we face the problem of finding some easy-to-verify conditions, which allow us to establish:
 \begin{enumerate}
 \item[-] $\mathfrak D$-stability of a given matrix $\mathbf A$;
 \item[-] the preservation of $\mathfrak D$-stability under some specific perturbations of ${\mathbf A}$, e.g. multiplication by a diagonal matrix.
 \end{enumerate}

Though an amount of research (see \cite{GUT2}, \cite{GUJU}) is due to ${\mathfrak D}$-stability (or eigenvalue clustering in a region ${\mathfrak D} \subset {\mathbb C}$), the progress in describing the matrix classes which are $\mathfrak D$-stable, and, moreover, preserve $\mathfrak D$-stability under multiplication by a specific diagonal matrix is still very little and concerns classical regions such as the left half-plane (see \cite{AM}, \cite{CROSS}, \cite{JOHN1}) and the unit disk (see \cite{KAB}). Here, we make certain efforts to fill this gap.

The paper is organized as follows. In Section 2, we study unbounded LMI regions and their properties, focusing on the most well-known examples, such as the shifted left half-plane, the conic sector, hyperbola and parabola. We also recall the concept of $({\mathfrak D},D)$-stability for an unbounded LMI region $\mathfrak D$. In Section 3, we introduce the crucial concept of this paper, namely, diagonal dominance with respect to a given LMI region $\mathfrak D$ (so-called diagonal $\mathfrak D$-dominance). We consider the particular cases of diagonal ${\mathfrak D}$-dominance, with respect to the most important LMI regions. Section 4 gives the basic results of the paper, i.e. the implications between diagonal ${\mathfrak D}$-dominance, ${\mathfrak D}$-stability and $({\mathfrak D},D)$-stability. Section 5 contains new results on the stability of the second-order systems, based on the Lyapunov theorem and the property of diagonal stability. Section 6 provides sufficient conditions for the eigenvalues of the matrix of a second-order system to be localized inside the shifted left half-plane (so-called ${\mathbb C}^-_{\alpha}$-stability). These conditions are based on the property of diagonal ${\mathfrak D}$-dominance. For some particular cases, we provide necessary and sufficient conditions of ${\mathbb C}^-_{\alpha}$-stability.  In Section 7, we study perturbations of second-order systems and introduce the concept of relative $D$-stability with the minimal decay rate $\alpha$. We provide sufficient conditions for $D$-stability and relative $D$-stability of second-order systems. In Section 8, we provide sufficient conditions for the stability of perturbed fractional-order systems. Finally, Section 9 contains numerical examples and calculations.

\subsection{Example} Here, we consider the example of the applications of results, basing on \cite{GUT2}, p. 23. Example 2.1. Let a feedback system (see \cite{GUT2}, p. 23, Figure 2.1) be described by the following open-loop transfer function:
$$G(s) = \frac{\omega_n^2}{s(s+2\zeta\omega_n)}, $$
where $\omega_n$ is the undamped natural frequency, $\zeta$ is the damping ratio, with the parameter-dependent feedback function of the form
$$H(s) = k_1s + k_2, $$
where $k_1$ and $k_2$ are the unknown positive parameters. We need to select the damping ratio $\zeta$ such that $\int_0^{\infty} y^2(t)dt$ is minimized for all positive values $k_1$ and $k_2$. We also require the asymptotic stability of the system.

The closed-loop transfer function of the system is
$$\frac{C(s)}{R(s)} = \frac{G(s)}{1+G(s)H(s)} = \frac{\omega_n^2}{s(s+2\zeta\omega_n) + \omega_n^2(k_1s + k_2)} =$$
$$\frac{\omega_n^2}{s^2+2\zeta\omega_n (1 + \frac{k_1 \omega_n}{2\zeta})s + \omega_n^2 k_2} = \ldots $$
Denoting $\hat{k}_1 := 1 + \frac{k_1 \omega_n}{2\zeta}$ and $\hat{k}_2 = k_2$, we obtain:
$$\ldots = \frac{\omega_n^2}{s^2+2\zeta\omega_n \hat{k}_1 s + \omega_n^2 \hat{k}_2}. $$

The state-space model of the closed-loop system is:
$$\dot{x} = \begin{pmatrix}-2\zeta\omega_n \hat{k}_1 & - \omega_n^2 \hat{k}_2 \\ 1 & 0\end{pmatrix}x; $$
$$ y = \begin{pmatrix}1 & 0 \end{pmatrix}x.$$

Note, that the parameter-dependent system matrix $$\widetilde{{\mathbf A}}(k):=\begin{pmatrix}-2\zeta\omega_n \hat{k}_1 & - \omega_n^2 \hat{k}_2 \\ 1 & 0\end{pmatrix}$$
is equivalent to the product of the form ${\mathbf D}\widetilde{{\mathbf A}}$, where
$${\mathbf D} = \begin{pmatrix}\hat{k}_1 & 0 \\ 0 & \frac{\hat{k}_2}{\hat{k}_1}\end{pmatrix}; \qquad \widetilde{{\mathbf A}} := \begin{pmatrix}-2\zeta\omega_n  & - \omega_n^2 \\ 1 & 0\end{pmatrix}. $$

Following the reasoning from \cite{GUT2}, we obtain for $f(\zeta) = \int_0^{\infty} y^2(t)dt$:
$$f(\zeta) = {\rm Tr}{\mathbf P}, $$
where ${\mathbf P}$ is chosen to be the solution of the Lyapunov equation:
$${\mathbf P}\widetilde{{\mathbf A}}(k) + (\widetilde{{\mathbf A}}(k))^T{\mathbf P} = -{\mathbf Q}, $$
where ${\mathbf Q} = {\mathbf e}{\mathbf e}^T$, ${\mathbf e}^T = (0, 1)$.
Simple calculations show that
$${\mathbf P} = \frac{1}{\hat{k}_2\omega_n}\begin{pmatrix}\frac{1}{\hat{k}_14\zeta\omega_n^2} & \frac{1}{ 2\omega_n} \\
\frac{1}{ 2\omega_n} & \hat{k}_1\zeta + \frac{\hat{k}_2}{\hat{k}_1 4\zeta} \end{pmatrix},$$
thus $$f(\zeta) =\frac{1}{\hat{k}_2\omega_n}\left(\frac{1}{\hat{k}_14\zeta\omega_n^2} + \hat{k}_1\zeta + \frac{\hat{k}_2}{\hat{k}_1 4\zeta}\right) = \frac{1}{\hat{k}_2\omega_n}\left( \hat{k}_1\zeta + \frac{1}{\hat{k}_14\zeta}\left(\hat{k}_2 + \frac{1}{\omega_n^2}\right)\right).$$
To minimize $f(\zeta)$, we solve the equation $f'(\zeta)=0$ and obtain
\begin{equation}\label{ratio}\zeta =  \frac{1}{2\hat{k}_1}\sqrt{\hat{k}_2 + \frac{1}{\omega_n^2}}.\end{equation} For asymptotic stability, we also require \begin{equation}\label{asst}\zeta\omega_n > 0.\end{equation}

The closed-loop poles of the system are $\lambda_{1,2} = x \pm iy,$ where
$$x = -\zeta\omega_n\hat{k}_1; \qquad y = \omega_n\sqrt{\hat{k}_2 - \zeta^2\hat{k}_1^2}.$$

Thus, for all positive values of parameters $\hat{k}_1, \hat{k}_2$, Equation \eqref{ratio} implies:
$$4\hat{k}_1^2\omega_n^2\zeta^2 = \hat{k}_2\omega_n^2 + 1. $$
Taking into account that  $$x^2 = \zeta^2\omega_n^2\hat{k}_1^2; \qquad y^2 = \omega_n^2\hat{k}_2 - x^2$$ and conditions of asymptotic stability \eqref{asst}, we obtain the left branch of the hyperbola
\begin{equation}\label{d}3x^2 - y^2 = 1. \end{equation}

Consider the region $\mathfrak D$, defined as follows:
$${\mathfrak D}: = H(a,b) = \{z = x+iy \in {\mathbb C}: x < - \frac{1}{|a|}; \ 1 - a^2x^2 + b^2y^2 < 0\}. $$

For some practical reasons, we may require the eigenvalues of the parameter-dependent system matrix $\widetilde{{\mathbf A}}(k)$ to be clustered in $\mathfrak D$ for all values of parameters $\hat{k}_1, \hat{k}_2$, satisfying $\hat{k}_2 > \hat{k}_1 > 1$. This leads us to a more general problem of $({\mathfrak D}, D)$-stability of $\widetilde{\mathbf A}$ (see Definition 5') with respect to the hyperbolic region $H(a,b)$, described in Example 3. Note, that the region ${\mathfrak D}$ bounded by the left brunch of the hyperbola \eqref{d} is exactly $H(3,1)$. As it is proved in Theorem \ref{mmain}, diagonal $H(a,b)$-dominance (see Definition 8) gives us sufficient conditions of localization of the spectrum of $\widetilde{{\mathbf A}}$ in the region $H(a,b)$. Then, diagonal ${\mathbb C}^{x_{max}}_{\theta_0}$-dominance (see Definition 7) with $x_{max} = -\frac{1}{|a|}$, $\theta_0 = \arctan|\frac{a}{b}|$ implies diagonal $H(a,b)$-dominance (by Lemma \ref{hyp}) and gives us sufficient condition for $(H(a,b), D)$-stability (by Theorem \ref{mainth}).

\section{Generalized $D$-stability in unbounded LMI regions}

\subsection{Definition and examples of LMI regions}
Consider the following type of regions, introduced in \cite{CHG} (see also \cite{CGA}).

{\bf Definition 4.} A subset ${\mathfrak D} \subset {\mathbb C}$ that can be defined as
\begin{equation}\label{LMI} {\mathfrak D} = \{z \in {\mathbb C}: \ {\mathbf L} + {\mathbf M}z+{\mathbf M}^T\overline{z} \prec 0\},\end{equation}
where ${\mathbf L}, {\mathbf M} \in {\mathcal M}^{n \times n}$, ${\mathbf L}^T = {\mathbf L}$, is called an {\it LMI region} with the {\it characteristic function} $f_{\mathfrak D}(z) = {\mathbf L} + z{\mathbf M}+\overline{z}{\mathbf M}^T$ and {\it generating matrices} $\mathbf M$ and $\mathbf L$. It was shown in \cite{KU}, that {\it an LMI region ${\mathfrak D}$ is open}. In the sequel, we shall use the notation $\overline{{\mathfrak D}}$ for its closure and $\partial{\mathfrak D}$ for its boundary.

Later on, we shall be especially interested in the following types of LMI regions.

{\bf Example 1.} The shifted left half-plane with the shift parameter $\alpha \in {\mathbb R}$
 $${\mathbb C}^-_{\alpha} = \{\lambda \in {\mathbb C} : {\rm Re}(\lambda) < \alpha\},$$
 with the characteristic function $$f_{{\mathbb C}^-_{\alpha}}(z) = z + \overline{z}-2\alpha.$$

{\bf Example 2.} The shifted conic sector around the negative direction of the real axis with the apex at the point $\alpha \in {\mathbb R}$ and the inner angle $2\theta$, $0 < \theta < \frac{\pi}{2}$,
$$
{\mathbb C}^{\alpha}_{\theta} = \{z = x+iy \in {\mathbb C}: (\alpha-x)\tan\theta > y > (-\alpha+x)\tan\theta\},
$$
with the characteristic function
$$f_{{\mathbb C}^{\alpha}_{\theta}} = \begin{pmatrix} \sin(\theta) & \cos(\theta) \\ - \cos(\theta) & \sin(\theta) \\ \end{pmatrix}z + \begin{pmatrix} \sin(\theta) & -\cos(\theta) \\  \cos(\theta) & \sin(\theta) \\ \end{pmatrix}\overline{z} - 2\alpha\sin(\theta)\begin{pmatrix} 1 & 0 \\  0 & 1 \\ \end{pmatrix}.$$

According to the above notation, the conic sector with the apex at the origin ($\alpha = 0$) will be denoted by ${\mathbb C}^{0}_{\theta}$.
In the sequel, we shall also consider the conic sector around the {\it positive} direction of the real axis with the apex at the origin and the inner angle $2\theta$, $0 < \theta < \frac{\pi}{2}$, denoting it as ${\mathbb C}^{+}_{\theta}$. Obviously, ${\mathbb C}^{+}_{\theta} = - {\mathbb C}^{0}_{\theta}$.

{\bf Example 3.} The left side of the left branch of the hyperbola
$$H(a,b) = \{z = x+iy \in {\mathbb C}: x < - \frac{1}{|a|}; \ 1 - a^2x^2 + b^2y^2 < 0\}, $$
with the characteristic function
$$f_{H(a,b)} = \begin{pmatrix} a & b \\ - b & a \\ \end{pmatrix}z + \begin{pmatrix} a & -b \\  b & a \\ \end{pmatrix}\overline{z} + 2\begin{pmatrix} 1 & 0 \\  0 & -1 \\ \end{pmatrix}.$$

{\bf Example 4.} The left side of the stability parabola $y^2 = -\epsilon^2 x$, where $\epsilon$ is a damping parameter
 $$P(\epsilon) = \{z = x+iy \in {\mathbb C}: \ y^2 < -\epsilon^2 x \}, $$
with the characteristic function
$$f_{P(\epsilon)} = \begin{pmatrix} \frac{1}{2} & -1 \\  0 & \frac{1}{2} \\ \end{pmatrix}z + \begin{pmatrix} \frac{1}{2} & 0 \\  -1 & \frac{1}{2} \\ \end{pmatrix}\overline{z} +  \begin{pmatrix} -\epsilon^2 & 0 \\  0 & 0 \\ \end{pmatrix}. $$

Let us recall some definitions and results from \cite{KU}. A cone ${\mathfrak D}_{rc}$ is called {\it the recession cone} of an unbounded LMI region $\mathfrak D$ if, for every $z \in {\mathfrak D}$, the set $z + {\mathfrak D}_{rc} \subseteq {\mathfrak D}$. Consider an unbounded LMI region $\mathfrak D$ with a generating matrix $\mathbf M$ in Formula \eqref{LMI} being nonnegative definite. Then the recession cone ${\mathfrak D}_{rc}$ of $\mathfrak D$ is known to be a closed convex conic sector around the negative direction of the real axis, with the apex at the origin. In the case when $\mathbf M$ is positive definite as it was shown in \cite{KU} (see \cite{KU}, Theorems 11 and 13), the interior of ${\mathfrak D}_{rc}$ is nonempty, i.e. the inner angle $\theta_0$ of ${\mathfrak D}_{rc}$ satisfies $0 < \theta_0 \leq \frac{\pi}{2}$. The inner angle $\theta_0$ of ${\mathfrak D}_{rc}$ can be calculated using the spectral characteristics of the generating matrices $\mathbf L$ and $\mathbf M$ (see \cite{KU}, Theorem 23). In the case when $\mathbf M$ is non-symmetric and its symmetric part $\frac{{\mathbf M} + {\mathbf M}^T}{2}$ is singular nonnegative definite, ${\mathfrak D}_{rc} = {\mathbb R}_{-}$, i.e. coincides with the negative direction of the real axis, including zero (see \cite{KU}, Corollary 8).

In the sequel, we shall use the following characteristics of an unbounded LMI region $\mathfrak D$:
\begin{enumerate}
\item[\rm 1.] The inner angle $\theta_0$ of the conic sector ${\mathfrak D}_{rc}$;
\item[\rm 2.] The value $x_{max} = \max\{x: x\in {\mathbb R}\cap {\mathfrak D}\}$. Note that $x_{max} \in \partial(\mathfrak D)$.
\end{enumerate}

Now let us calculate $x_{max}$, ${\mathfrak D}_{rc}$ and its inner angle $\theta_0$ for the regions, considered in Examples 1-4.

\begin{enumerate}
\item[\rm 1.] For ${\mathfrak D} = {\mathbb C}^-_{\alpha}$, it follows directly from the definitions that $x_{max} = \alpha$, ${\mathfrak D}_{rc} = \overline{{\mathbb C}_0^-}$ and $\theta_0 = \frac{\pi}{2}$.
\item[\rm 2.] For ${\mathfrak D} = {\mathbb C}^{\alpha}_{\theta}$, it is easy to see that $x_{max} = \alpha$, ${\mathfrak D}_{rc} = \overline{{\mathbb C}^{0}_{\theta}}$ and $\theta_0 = \theta$.
\item[\rm 3.] For ${\mathfrak D} = H(a,b)$, using the results from \cite{KU}, we obtain that $x_{max} = -\frac{1}{|a|}$, ${\mathfrak D}_{rc} = \overline{{\mathbb C}^{0}_{\theta_0}}$, where $\theta_0 = \arctan|\frac{a}{b}|$.
\item[\rm 4.] For ${\mathfrak D} = P(\epsilon)$, it was shown in \cite{KU} (see \cite{KU}, Section 7, Subsection 7.7), that $x_{max} = 0$,  ${\mathfrak D}_{rc} = {\mathbb R}_-$, i.e. the negative direction of the real axis.
\end{enumerate}

The calculation of $x_{max}$ of an arbitrary LMI region $\mathfrak D$ in terms of the spectral characteristics of the generating matrices $\mathbf L$ and $\mathbf M$ is given in \cite{KU}, Theorem 27.

In the sequel, in some cases we shall replace the study of an arbitrary unbounded LMI region $\mathfrak D$ with the study of the (open) shifted conic sector \begin{equation}\label{sect}{\mathbb C}^{x_{max}}_{\theta_0}= {\rm int}(x_{max} + {\mathfrak D}_{rc}) \subseteq {\mathfrak D}. \end{equation}
Inclusion \eqref{sect} follows from the inclusion $x_{max} \in \overline{{\mathfrak D}}$ and the fact that the recession cone ${\mathfrak D}_{rc}$ of an LMI region $\mathfrak D$ coincides with the recession cone of $\overline{\mathfrak D}$ (see \cite{KU}, Corollary 6).

\subsection{$({\mathfrak D}, D)$-stability of matrices}
In \cite{KUPA}, the following generalizations of the concept of $D$-stability to the case of an unbounded LMI region $\mathfrak D$ were considered (separately for the cases when $0 \in \overline{{\mathfrak D}}$ and when $0 \in {\mathbb C}\setminus \overline{{\mathfrak D}}$).

 {\bf Definition 5.} Given an unbounded LMI region $\mathfrak D$ with $0 \in \overline{{\mathfrak D}}$, we say that an $n \times n$ real matrix $\mathbf A$ is {\it (multiplicative) $D$-stable with respect to $\mathfrak D$} or simply {\it (multiplicative) $({\mathfrak D}, D)$-stable} if $\sigma({\mathbf D}{\mathbf A}) \subset {\mathfrak D}$ for every ${\mathbf D} \in {\mathcal D}^+$.

 {\bf Definition 5'.} Given an unbounded LMI region $\mathfrak D$ with $0 \in {\mathbb C} \setminus \overline{{\mathfrak D}}$, we say that an $n \times n$ real matrix $\mathbf A$ is {\it (multiplicative) $D$-stable with respect to $\mathfrak D$} or simply {\it (multiplicative) $({\mathfrak D}, D)$-stable} if $\sigma({\mathbf D}{\mathbf A}) \subset {\mathfrak D}$ for every ${\mathbf D} \in {\mathcal D}^+_{\geq 1}$.

Now let us consider the following particular cases of $({\mathfrak D}, D)$-stability.

\begin{enumerate}
\item[\rm 1.]{\bf $({\mathbb C}^-_{\alpha}, D)$-stability.} Here, we consider separately the case, when $\alpha \geq 0$ which implies $0 \in \overline{{\mathbb C}^-_{\alpha}}$ and $\alpha < 0$ which implies $0 \in {\mathbb C} \setminus \overline{{\mathbb C}^-_{\alpha}}$. We say that an $n \times n$ real matrix $\mathbf A$ is {\it (multiplicative) $({\mathbb C}^-_{\alpha}, D)$-stable} if ${\rm Re}(\lambda )<\alpha$ for all $\lambda \in
\sigma (\mathbf{DA})$, where ${\mathbf D}$ is any matrix from ${\mathcal D}^+$ (if $\alpha \geq 0$) or ${\mathcal D}^+_{\geq 1}$ (if $\alpha < 0$).
\item[\rm 2.] {\bf $({\mathbb C}^0_{\theta}, D)$-stability.} For this case, we have $0 \in \overline{{\mathbb C}^0_{\theta}}$. For a given value $\theta$, $0 < \theta < \frac{\pi}{2}$, we call an $n \times n$ real matrix $\mathbf A$ {\it (multiplicative) $({\mathbb C}^0_{\theta}, D)$-stable} if $\sigma({\mathbf D}{\mathbf A}) \subset {\mathbb C}^0_{\theta}$ for every positive diagonal matrix $\mathbf D$.
\item[\rm 3.]{\bf $(H(a,b), D)$-stability.}  For this case, we have $0 \in {\mathbb C}\setminus \overline{H(a,b)}$. We say that an $n \times n$ real matrix $\mathbf A$ is {\it (multiplicative) $(H(a,b), D)$-stable} if $\sigma({\mathbf D}{\mathbf A}) \subset H(a,b)$ for every matrix ${\mathbf D} \in {\mathcal D}^+_{\geq 1}$.
\item[\rm 4.]{\bf $(P(\epsilon), D)$-stability.} Here, $0 \in \overline{P(\epsilon)}$, and we say that $n \times n$ real matrix $\mathbf A$ is {\it (multiplicative) $(P(\epsilon), D)$-stable} if $\sigma({\mathbf D}{\mathbf A}) \subset P(\epsilon)$ for every positive diagonal matrix ${\mathbf D}$.
\end{enumerate}
\section{Extended diagonal dominance with respect to an LMI region $\mathfrak D$}

\subsection{General definition}
 Given a (bounded or unbounded) LMI region ${\mathfrak D} \subset {\mathbb C}$ (e.g. an LMI region). Let the intersection of ${\mathfrak D}$ with the real axis $\mathbb R$ be denoted by
 $$(\alpha,\beta): = {\mathfrak D}\cap{\mathbb R},$$
  where both $\alpha$ and $\beta$ may be infinite. According to this notation, $\beta = x_{max}$. Note, that ${\mathfrak D}\cap{\mathbb R}$ is nonempty whenever $\mathfrak D$ is nonempty (see \cite{KU}, Lemma 21).

Let the part of the boundary of $\mathfrak D$ which lies above the real axis (if $\mathfrak D$ is bounded from above) be represented by a concave function $y(x): (\alpha, \beta)\rightarrow {\mathbb R}_{+}$. Note, that the boundary function $y(x)$ of an LMI region $\mathfrak D$, is defined implicitly.

 We define the function $r(x): (\alpha,\beta)\rightarrow {\mathbb R}^+$ as follows:
 \begin{enumerate}
\item[\rm Case I.] Let either $\alpha$ or $\beta$ be finite, $\mathfrak D$ be bounded from above. Then
$$r(x) := \frac{\min(|x-\alpha|,|x-\beta|)\cdot y(x)}{\sqrt{(\min(|x-\alpha|,|x-\beta|))^2 + y^2(x)}}. $$
\item[\rm Case II.] Let both $\alpha$ and $\beta$ be infinite (in this case, $(\alpha,\beta) = {\mathbb R}$ and $\mathfrak D$ is a horizontal stripe, see \cite{KU}, Theorem 12, Part 2). Then
 $$r(x) : = \lim_{\alpha,\beta\rightarrow \infty}\frac{\min(|x-\alpha|,|x-\beta|)\cdot y(x)}{\sqrt{(\min(|x-\alpha|,|x-\beta|))^2 + y^2(x)}} = y(x) = {\rm const}. $$
\item[\rm Case III.] Let $\mathfrak D$ be unbounded from above (in this case, $\mathfrak D$ is a vertical stripe or a halfplane, see \cite{KU}, Theorem 12, Part 1). Then
$$r(x) := \lim_{y\rightarrow \infty}\frac{\min(|x-\alpha|,|x-\beta|)\cdot y}{\sqrt{(\min(|x-\alpha|,|x-\beta|))^2 + y^2}} = \min(|x-\alpha|,|x-\beta|). $$
\end{enumerate}

{\bf General definition.} A matrix ${\mathbf A} \in {\mathcal M}^{n \times n}$ is called {\it strictly row diagonally dominant with respect to the region $\mathfrak D$} or simply {\it diagonally ${\mathfrak D}$-dominant} if it satisfies the following properties:
\begin{enumerate}
\item[\rm 1.] $a_{ii} \in (\alpha,\beta)$ for all $i = 1, \ \ldots, \ n$;
\item[\rm 2.] $r(a_{ii}) > \sum_{i \neq j}|a_{ij}|$ for all $i = 1, \ \ldots, \ n$.
\end{enumerate}
A matrix ${\mathbf A} \in {\mathcal M}^{n \times n}$ is called {\it strictly column diagonally dominant with respect to the region $\mathfrak D$} if ${\mathbf A}^T$ is diagonally ${\mathfrak D}$-dominant.

Note, that the above definition is meaningful: given an arbitrary LMI region ${\mathfrak D} \subset {\mathbb C}$, for example, any diagonal matrix with principal diagonal entries from $(\alpha,\beta) = {\mathfrak D}\cap{\mathbb R}$, will be necessarily diagonally ${\mathfrak D}$-dominant.

\subsection{Particular cases}
First, consider the shifted half-plane ${\mathbb C}^-_{\alpha}$, with the shift parameter $\alpha \in {\mathbb R}$.

 {\bf Definition 6.} Given a value $\alpha \in {\mathbb R}$, a real $n \times n$ matrix ${\mathbf A}$ is called {\it diagonally ${\mathbb C}^-_{\alpha}$-dominant} if the following inequalities hold:
\begin{enumerate}
\item[\rm 1.]$|a_{ii} - \alpha| > \sum_{j\neq i}|a_{ij}| \qquad i = 1, \ \ldots, \ n.$
\item[\rm 2.] $a_{ii} < \alpha$, $i = 1, \ \ldots, \ n,$
\end{enumerate}

For the case $\alpha = 0$, the Definition 6 gives the class of diagonally dominant matrices (see Definition 1) with negative principal diagonal entries, i.e. a subclass of NDD matrices.

The second region we consider, is the shifted conic sector around the negative direction of the real axis with the apex at the point $\alpha \in {\mathbb R}$ and the inner angle $2\theta$, $0 < \theta < \frac{\pi}{2}$. In this case, it easily follows from geometrical reasoning that $r(x) = \sin \theta|x - \alpha|$.

{\bf Definition 7.} Given two values $\theta \in (0, \frac{\pi}{2}]$ and $\alpha \in {\mathbb R}$, a matrix ${\mathbf A} \in {\mathcal M}^{n \times n}$ is called {\it diagonally ${\mathbb C}^{\alpha}_{\theta}$-dominant} if the following inequalities hold:
\begin{enumerate}
\item[\rm 1.]$\sin \theta|a_{ii} - \alpha| > \sum_{j\neq i}|a_{ij}| \qquad i = 1, \ \ldots, \ n.$
\item[\rm 2.] $a_{ii} < \alpha$, $i = 1, \ \ldots, \ n.$
\end{enumerate}

For the conic sector ${\mathbb C}^{0}_{\theta}$, with the apex at the origin, we have the following definition.

{\bf Definition 7'.} Given a value $\theta \in (0, \frac{\pi}{2}]$, a real $n \times n$ matrix ${\mathbf A}$ is called {\it diagonally ${\mathbb C}^{0}_{\theta}$-dominant} if the following inequalities hold:
 \begin{enumerate}
\item[\rm 1.]$\sin \theta|a_{ii}| > \sum_{j\neq i}|a_{ij}| \qquad i = 1, \ \ldots, \ n.$
\item[\rm 2.] $a_{ii} < 0$, $i = 1, \ \ldots, \ n.$
\end{enumerate}

The following extension of Definition 7' was provided in \cite{NIS2} (see also \cite{NIS}).

{\bf Definition 7''.} A matrix ${\mathbf A} \in {\mathcal M}^{n \times n}({\mathbb C})$ is called {\it generalized diagonally dominant with a given strength factor $\tau$}, where $0 < \tau < 1$, if there exist positive scalars (weights) $m_i$, $i = 1, \ \ldots, \ n$ such that
$$\tau m_i|a_{ii}|>\sum_{j\neq i }m_j|a_{ij}|, \qquad i = 1, \ \ldots, \ n. $$
If, in addition, $a_{ii} < 0$, $i = 1, \ \ldots, \ n,$ then ${\mathbf A}$ is called {\it negative diagonally dominant (NDD) with strength factor $\tau$}. Here, we may put $\sin \theta:=\tau$ and obtain that ${\mathbf A}$ is NDD with strength factor $\tau$ if and only if there is a positive diagonal matrix ${\mathbf M}$ such that ${\mathbf A}{\mathbf M}$ is diagonally ${\mathbb C}^{0}_{\theta}$-dominant.

The following definition corresponds to the hyperbolic region $H(a,b)$. Here, we have $\alpha = -\frac{1}{|a|}$, $y(x) = \frac{\sqrt{a^2x^2 - 1}}{|b|}$, thus
$$r(x) = \frac{|x+\frac{1}{|a|}|\sqrt{a^2x^2 - 1}}{\sqrt{b^2(x+\frac{1}{|a|})^2 + (a^2x^2-1)}} =\frac{(-|a|x-1)\sqrt{1 -|a|x)}}{\sqrt{b^2(-|a|x-1)+a^2(1-|a|x)}}. $$

{\bf Definition 8.} Given two values $a,b \in {\mathbb R}\setminus\{0\}$, a real $n \times n$ matrix ${\mathbf A}$ is called {\it diagonally $H(a,b)$-dominant} if the following inequalities hold:
 \begin{enumerate}
\item[\rm 1.]$\frac{(-|a|a_{ii}-1)\sqrt{1 -|a|a_{ii})}}{\sqrt{b^2(-|a|a_{ii}-1)+a^2(1-|a|a_{ii})}} > \sum_{j\neq i}|a_{ij}| \qquad i = 1, \ \ldots, \ n.$
\item[\rm 2.] $a_{ii} < -\frac{1}{|a|}$, $i = 1, \ \ldots, \ n.$
\end{enumerate}

Finally, the following definition corresponds to the stability parabola $P(\epsilon)$. In this case,
$$r(x) = \frac{|x||\epsilon|\sqrt{-x}}{\sqrt{x^2 - \epsilon^2x}}= \frac{|\epsilon x|}{\sqrt{\epsilon^2 -x}} $$

{\bf Definition 9.} Given a value $\epsilon \neq 0$, a real $n \times n$ matrix ${\mathbf A}$ is called {\it diagonally $P(\epsilon)$-dominant} if the following inequalities hold:
 \begin{enumerate}
\item[\rm 1.]$\frac{|\epsilon a_{ii}|}{\sqrt{\epsilon^2 -a_{ii}}} > \sum_{j\neq i}|a_{ij}| \qquad i = 1, \ \ldots, \ n.$
\item[\rm 2.] $a_{ii} < 0$, $i = 1, \ \ldots, \ n.$
\end{enumerate}

\subsection{Properties of diagonally $\mathfrak D$-dominant matrices}
Here, we consider the basic properties of diagonally $\mathfrak D$-dominant matrices, that we shall use later. First, let us state and prove the following technical lemma.
\begin{lemma}\label{vosp}
Let the function $r: ({\mathbb R}_+ \times {\mathbb R}_+) \rightarrow {\mathbb R}$ be defined by the formula:
$$r(a,b):=\frac{ab}{\sqrt{a^2+b^2}}. $$
Then $r(a_2, b_2) \leq r(a_1, b_1)$ whenever $a_2 \leq a_1$, $b_2 \leq b_1$.
\end{lemma}
{\bf Proof.} Since $a_2 \leq a_1 $ and $b_2 \leq b_1$, we obtain the following inequalities:
$$a^2_1a^2_2(b_2^2 - b_1^2) + b^2_1b^2_2(a_2^2 - a_1^2) \leq 0;$$
$$a^2_2b_2^2(a^2_1 + b_1^2) - a^2_1b^2_1(a_2^2 + b_2^2) \leq 0;$$
$$a^2_2b_2^2(a^2_1 + b_1^2) \leq a^2_1b^2_1(a_2^2 + b_2^2);$$
$$\frac{a^2_2b_2^2}{a^2_2 + b_2^2} \leq \frac{a^2_1b^2_1}{a_1^2 + b_1^2};$$
$$\frac{a_2b_2}{\sqrt{a^2_2 + b_2^2}} \leq \frac{a_1b_1}{\sqrt{a_1 + b_1}}.$$
The last inequality means $r(a_2, b_2) \leq r(a_1, b_1)$.
$\square$

The geometrical meaning of $r(a,b)$ is as follows: it provides the length of the altitude to the hypothenuse in the right triangle with the lengths of the other two sides to be $a$ and $b$.

\begin{lemma}\label{inclu} Let ${\mathfrak D}_1, {\mathfrak D}_2 \subset {\mathbb C}$ be two nonempty LMI regions satisfying the inclusion ${\mathfrak D}_2 \subseteq {\mathfrak D}_1$. Then any diagonally ${\mathfrak D}_2$-dominant matrix is necessarily diagonally ${\mathfrak D}_1$-dominant. In particular, all diagonally ${\mathbb C}^-_{\alpha}$-dominant matrices with $\alpha \leq 0$ and all diagonally ${\mathbb C}^0_{\theta}$-dominant matrices are diagonally dominant.
\end{lemma}
{\bf Proof.} Let ${\mathfrak D}_2 \subset {\mathfrak D}_1$. Since ${\mathfrak D}_1$, ${\mathfrak D}_2$ are LMI regions, they are both convex and symmetric with respect to the real axis. Thus, by \cite{KU}, Lemma 21, ${\mathfrak D}_1 \cap {\mathbb R} = (\alpha_1, \beta_1) \neq \emptyset$ and ${\mathfrak D}_2 \cap {\mathbb R} = (\alpha_2, \beta_2) \neq \emptyset$. Since ${\mathfrak D}_2 \subset {\mathfrak D}_1$ we conclude $(\alpha_2, \beta_2) \subseteq (\alpha_1, \beta_1)$, i.e. $\alpha_2 \geq \alpha_1$ and $\beta_2 \leq \beta_1$. Let $\mathbf A$ be a diagonally ${\mathfrak D}_2$-dominant matrix. Then it satisfies the following properties:
\begin{enumerate}
\item[\rm 1.] $a_{ii} \in (\alpha_2,\beta_2)$ for all $i = 1, \ \ldots, \ n$;
\item[\rm 2.] $r_2(a_{ii}) > \sum_{i \neq j}|a_{ij}|$ for all $i = 1, \ \ldots, \ n$.
\end{enumerate}
First, the inclusions $a_{ii} \in (\alpha_2,\beta_2)$ and $(\alpha_2, \beta_2) \subseteq (\alpha_1, \beta_1)$ imply $a_{ii} \in (\alpha_1,\beta_1)$. Then we need to show that $r_1(a_{ii}) \geq r_2(a_{ii}) > \sum_{i \neq j}|a_{ij}|$.
We consider the following cases.
\begin{enumerate}
\item[\rm Case I.] Let either $\alpha_2$ or $\beta_2$ be finite, ${\mathfrak D}_2$ be bounded from above. Then
$$r_2(a_{ii}) = \frac{\min(|a_{ii}-\alpha_2|,|a_{ii}-\beta_2|)\cdot y_2(a_{ii})}{\sqrt{(\min(|a_{ii}-\alpha_2|,|a_{ii}-\beta_2|))^2 + y_2^2(a_{ii})}}. $$
If either $\alpha_1$ or $\beta_1$ is finite and ${\mathfrak D}_1$ is bounded from above, we have $$r_1(a_{ii}) = \frac{\min(|a_{ii}-\alpha_1|,|a_{ii}-\beta_1|)\cdot y_1(a_{ii})}{\sqrt{(\min(|a_{ii}-\alpha_1|,|a_{ii}-\beta_1|))^2 + y_1^2(a_{ii})}}. $$
Then, ${\mathfrak D}_2 \subset {\mathfrak D}_1$ obviously implies $y_2(a_{ii}) \leq y_1(a_{ii})$ for every $a_{ii} \in (\alpha_2,\beta_2)$. The inequalities $\alpha_2 \geq \alpha_1$ and $\beta_2 \leq \beta_1$ imply $\min(|a_{ii}-\alpha_2|,|a_{ii}-\beta_2|) \leq \min(|a_{ii}-\alpha_1|,|a_{ii}-\beta_1|)$ for every $a_{ii} \in (\alpha_2,\beta_2)$. Applying Lemma \ref{vosp}, we obtain $r_2(a_{ii}) \leq r_1(a_{ii})$ for every $a_{ii} \in (\alpha_2,\beta_2)$. Now let us consider the case when $(\alpha_1, \beta_1) = {\mathbb R}$, i.e. ${\mathfrak D}_1$ is a horizontal stripe. In this case, $r_1(a_{ii}) = y_1(a_{ii}) = {\rm const}$, and we obtain:
$$r_2(a_{ii}) = \frac{\min(|a_{ii}-\alpha_2|,|a_{ii}-\beta_2|)\cdot y_2(a_{ii})}{\sqrt{(\min(|a_{ii}-\alpha_2|,|a_{ii}-\beta_2|))^2 + y_2^2(a_{ii})}} \leq y_2(a_{ii}) \leq y_1(a_{ii}) = r_1(a_{ii}).$$
Finally, let ${\mathfrak D}_1$ be unbounded from above. Then we have $r_1(a_{ii}) = \min(|a_{ii}-\alpha|,|a_{ii}-\beta|)$, and
$$r_2(a_{ii}) = \frac{\min(|a_{ii}-\alpha_2|,|a_{ii}-\beta_2|)\cdot y_2(a_{ii})}{\sqrt{(\min(|a_{ii}-\alpha_2|,|a_{ii}-\beta_2|))^2 + y_2^2(a_{ii})}} \leq \min(|a_{ii}-\alpha_2|,|a_{ii}-\beta_2|) $$ $$ \leq \min(|a_{ii}-\alpha_1|,|a_{ii}-\beta_1|) = r_1(a_{ii}).$$
\item[\rm Case II.] Let both $\alpha_2$ and $\beta_2$ be infinite (${\mathfrak D}_2$ is a horizontal stripe). Since ${\mathfrak D}_2 \subseteq {\mathfrak D}_1$, it is easy to see that ${\mathfrak D}_1$ is also a horizontal stripe.
     Then
 $$r_2(a_{ii}) = y_2(a_{ii})\leq y_1(a_{ii}) = r_1(a_{ii}). $$
\item[\rm Case III.] The regions ${\mathfrak D}_2$ is unbounded from above. In this case, the inclusion ${\mathfrak D}_2 \subseteq {\mathfrak D}_1$ implies that ${\mathfrak D}_1$ is also unbounded from above. We have
     $$r_2(a_{ii}) = \min(|a_{ii}-\alpha_2|,|a_{ii}-\beta_2|)\leq  \min(|a_{ii}-\alpha_1|,|a_{ii}-\beta_1|) = r_1(a_{ii}). $$
\end{enumerate}
$\square$

From Lemma \ref{inclu}, we deduce the following statement.
\begin{lemma} Given an unbounded LMI region ${\mathfrak D}$, let the conic sector ${\mathbb C}^{x_{max}}_{\theta_0}$, defined by Formula \eqref{sect}, be nonempty. Then any diagonally ${\mathbb C}^{x_{max}}_{\theta_0}$-dominant matrix is necessarily diagonally $\mathfrak D$-dominant. In particular, if $a,b \in {\mathbb R}\setminus \{0\}$, $x_{max} = -\frac{1}{|a|}$, $\theta_0 = \arctan|\frac{a}{b}|$, then any diagonally ${\mathbb C}^{x_{max}}_{\theta_0}$-dominant matrix is diagonally $H(a,b)$-dominant.
\end{lemma}
{\bf Proof.} The proof follows from the inclusion ${\mathbb C}^{x_{max}}_{\theta_0}\subseteq {\mathfrak D}$ (see Formula \eqref{sect}) and Lemma \ref{inclu}. $\square$

Now let us consider the following question: if $\mathbf A$ is diagonally $\mathfrak D$-dominant (with respect to some LMI region $\mathfrak D$), ${\mathbf D}$ is positive diagonal, when ${\mathbf D}{\mathbf A}$ be also diagonally $\mathfrak D$-dominant? The answer to this question will lead to sufficient conditions of $({\mathfrak D}, D)$-stability. The following result is obvious.

\begin{lemma}\label{obv} Let $\mathfrak D$ be an unbounded LMI region such that ${\mathfrak D}\cap{\mathbb R} = {\mathbb R}_-$ and the corresponding function $r(x)$ is uniform, i.e. $r(ax) = ar(x)$ for any $x < 0$, $a > 0$. If ${\mathbf A} \in {\mathcal M}^{n \times n}$ is diagonally $\mathfrak D$-dominant, then ${\mathbf D}{\mathbf A}$ is also diagonally $\mathfrak D$-dominant for any positive diagonal matrix $\mathbf D$. In particular, if $\mathbf A$ is diagonally ${\mathbb C}^{0}_{\theta}$-dominant for any $\theta \in (0, \frac{\pi}{2}]$, then ${\mathbf D}{\mathbf A}$ is also diagonally ${\mathbb C}^{0}_{\theta}$-dominant for any positive diagonal matrix $\mathbf D$.
\end{lemma}

Now let us consider more general cases when $r(x)$ is not uniform.

\begin{lemma}\label{par} Let $\mathfrak D$ be an unbounded LMI region such that ${\mathfrak D}\cap{\mathbb R} = {\mathbb R}_-$. If ${\mathbf A} \in {\mathcal M}^{n \times n}$ is diagonally $\mathfrak D$-dominant, then ${\mathbf D}{\mathbf A}$ is also diagonally $\mathfrak D$-dominant for any ${\mathbf D} \in {\mathcal D}^+_{(0,1]}$.
\end{lemma}

{\bf Proof.} Let ${\mathbf A}$ be diagonally $\mathfrak D$-dominant. Then it satisfies the following conditions:
 \begin{enumerate}
\item[\rm 1.] $a_{ii} \in (\alpha,\beta)$ for all $i = 1, \ \ldots, \ n$;
\item[\rm 2.] $r(a_{ii}) > \sum_{i \neq j}|a_{ij}|$ for all $i = 1, \ \ldots, \ n$.
\end{enumerate} Taking into account that ${\mathfrak D}\cap{\mathbb R} = (\alpha,\beta) = (-\infty,0)$, we obtain $a_{ii} < 0$ and the corresponding function $r(a_{ii})$ be defined as follows:
$$r(a_{ii}) = \frac{-a_{ii} y(a_{ii})}{\sqrt{a_{ii}^2 + y^2(a_{ii})}}. $$
Since the function $y(x)$ is concave whenever $x < 0$, we easily obtain the inequality
\begin{equation}\label{conc}y(d_{ii}a_{ii}) \geq d_{ii}y(a_{ii}), \end{equation} whenever $d_{ii} \in [0,1]$, $i = 1, \ \ldots, \ n$.

 Then, checking the definition of diagonal $\mathfrak D$-dominance for the matrix ${\mathbf D}{\mathbf A}$, we obtain: $d_{ii}a_{ii} < 0$, whenever $d_{ii}>0$, $a_{ii} < 0$, $i = 1, \ \ldots, \ n$.  Applying Inequality \eqref{conc} and Lemma \ref{vosp}, we obtain the following estimates for $r(d_{ii}a_{ii})$, $d_{ii} \in (0,1]$:
 $$r(d_{ii}a_{ii}) = \frac{-d_{ii}a_{ii} y(d_{ii}a_{ii})}{\sqrt{(d_{ii}a_{ii})^2 + y^2(d_{ii}a_{ii})}} \geq  \frac{-d_{ii}a_{ii} d_{ii}y(a_{ii})}{\sqrt{d^2_{ii}a_{ii}^2 + d^2_{ii}y^2(a_{ii})}} = $$
 $$\frac{-d_{ii}a_{ii} y(a_{ii})}{\sqrt{a_{ii}^2 + y^2(a_{ii})}} = d_{ii}r(a_{ii}) > d_{ii}\sum_{i \neq j}|a_{ij}| = \sum_{i \neq j}|d_{ii}a_{ij}|,$$
 for all $i = 1, \ \ldots, \ n$. Thus ${\mathbf D}{\mathbf A}$ is diagonally $\mathfrak D$-dominant for any ${\mathbf D} \in {\mathcal D}^+_{(0,1]}$.
$\square$

Now let us consider the case of a shifted conic sector, which is of major practical importance.

\begin{lemma}\label{hyp} Given two values $\theta \in (0, \frac{\pi}{2}]$ and $\alpha < 0$. If ${\mathbf A} \in {\mathcal M}^{n \times n}$ is diagonally ${\mathbb C}^{\alpha}_{\theta}$-dominant, then ${\mathbf D}{\mathbf A}$ is also diagonally ${\mathbb C}^{\alpha}_{\theta}$-dominant for any ${\mathbf D} \in {\mathcal D}_{\geq 1}$.
\end{lemma}
{\bf Proof.} Since $\mathbf A$ is diagonally ${\mathbb C}^{\alpha}_{\theta}$-dominant, we have by Definition 7:
\begin{enumerate}
\item[\rm 1.]$\sin \theta|a_{ii} - \alpha| > \sum_{j\neq i}|a_{ij}| \qquad i = 1, \ \ldots, \ n.$
\item[\rm 2.] $a_{ii} < \alpha < 0$, $i = 1, \ \ldots, \ n.$
\end{enumerate}
Thus, for the matrix ${\mathbf D}{\mathbf A}$, we obtain: $d_{ii}a_{ii} < d_{ii}\alpha < \alpha < 0$, for every $d_{ii} \geq 1$, $i = 1, \ \ldots, \ n$. Then,
$$\sin \theta|d_{ii}a_{ii} - \alpha| = d_{ii}\sin \theta\left|\frac{\alpha}{d_{ii}} - a_{ii}\right| \geq d_{ii}\sin \theta|\alpha - a_{ii}| > d_{ii}\sum_{j\neq i}|a_{ij}| = \sum_{j\neq i}|d_{ii}a_{ij}|, $$
for every $d_{ii} \geq 1$, $i = 1, \ \ldots, \ n$. The above inequalities show that ${\mathbf D}{\mathbf A}$ is diagonally ${\mathbb C}^{\alpha}_{\theta}$-dominant for any ${\mathbf D} \in {\mathcal D}_{\geq 1}$.
 $\square$

\section{Diagonal $\mathfrak D$-dominance implies $\mathfrak D$-stability and $({\mathfrak D}, D)$-stability}

First, we recall the following fundamental statement of matrix eigenvalue localization (see, for example, \cite{HOJ}, p. 344).
\begin{theorem}[Gershgorin]\label{GER} Let ${\mathbf A} = \{a_{i,j}\}_{i,j = 1}^n \in {\mathcal M}^{n \times n}({\mathbb C})$, define
$$R_i := \sum_{j =1; \ i \neq j}^n |a_{ij}|, \qquad 1 \leq i \leq n.$$ Let $D(a_{ii},R_i) \subset {\mathbb C}$ be a closed disk centered at $a_{ii}$ with the radius $R_i$. Then all the eigenvalues of ${\mathbf A}$ are located in the union of $n$ discs $$G({\mathbf A}):= \bigcup_{i=1}^n D(a_{ii},R_i). $$
\end{theorem}
A disk $D(a_{ii},R_i)$ defined in the statement of Theorem \ref{GER}, is called a {\it Gershgorin disk}.

Let us recall a well-known fact that if a complex matrix ${\mathbf A} = \{a_{ij}\}_{i,j = 1}^n$ is strictly diagonally dominant and its principal diagonal entries satisfy $a_{ii}<0 , \ i= 1, \ \ldots, \ n,$ then $\mathbf A$ is Hurwitz stable, i.e. $\sigma({\mathbf A}) \subset {\mathbb C}^-$. For the proof, it is enough to notice that the union of all Gershgorin discs $G({\mathbf A})$, which contains $\sigma({\mathbf A})$, lies entirely in ${\mathbb C}^-$. To generalize this result, we first recall the following definition (see, for example, \cite{KU2}).

{\bf Definition 10.} Let a set ${\mathfrak D} \subset {\mathbb C}$ be symmetric with respect to the real axis. A matrix ${\mathbf A} \in {\mathcal M}^{n \times n}$ is called {\it stable with respect to $\mathfrak D$} or simply {\it $\mathfrak D$-stable} if $\sigma({\mathbf A}) \subset {\mathfrak D}$. In this case, ${\mathfrak D}$ is called a {\it stability region}.

Basing on Gershgorin's theorem, we obtain the following statement.

\begin{theorem}\label{mmain} Given a (bounded or unbounded) LMI region ${\mathfrak D} \subset {\mathbb C}$. If a matrix ${\mathbf A} \subset {\mathcal M}^{n \times n}$ is diagonally ${\mathfrak D}$-dominant then it is $\mathfrak D$-stable.
\end{theorem}
{\bf Proof.} Let us consider Gershgorin disks $D(a_{ii},R_i)$, $i = 1, \ \ldots, \ n$. For the proof, it will be enough to show the inclusions
\begin{equation}\label{i} D(a_{ii},R_i) \subseteq {\mathfrak D}, \qquad i = 1, \ \ldots, \ n,\end{equation}
and to apply Gershgorin's theorem. Since $\mathbf A$ is diagonally $\mathfrak D$-dominant, by definition we have:
$$a_{ii} \in (\alpha,\beta) = {\mathfrak D}\cap{\mathbb R} \qquad i = 1, \ \ldots, \ n, $$
i.e. the center of each Gershgorin disk $D(a_{ii},R_i)$ lies in $\mathfrak D$. Now let us estimate their radii $R_i = \sum_{j =1; \ i \neq j}^n |a_{ij}|,$ $1 \leq i \leq n.$ Again by definition of diagonal $\mathfrak D$-dominance, we have
$$R_i = \sum_{j =1; \ i \neq j}^n |a_{ij}| < r(a_{ii}), \qquad i = 1, \ \ldots, \ n. $$
Hence $D(a_{ii},R_i) \subset D(a_{ii},r(a_{ii}))$, where $D(a_{ii},r(a_{ii}))$ is a closed disk centered at $a_{ii}$ with the radius $r(a_{ii})$.
Obviously, for Inclusion \eqref{i} to hold, it is enough to show that
$$D(a_{ii}, r(a_{ii})) \subset \overline{{\mathfrak D}}, \qquad i = 1, \ \ldots, \ n, $$
where $\overline{{\mathfrak D}}$ is the closure of $\mathfrak D$.
For the proof, we consider the following cases.
\begin{enumerate}
\item[\rm Case I.] Let either $\alpha$ or $\beta$ be finite, $\mathfrak D$ be bounded from above. Without loss the generality of the reasoning, we may assume that both $\alpha$ and $\beta$ are finite. Then
$$r(a_{ii}) = \frac{\min(|a_{ii}-\alpha|,|a_{ii}-\beta|)\cdot y(a_{ii})}{\sqrt{(\min(|a_{ii}-\alpha|,|a_{ii}-\beta|))^2 + y^2(a_{ii})}}. $$
Let $\min(|a_{ii}-\alpha|,|a_{ii}-\beta|) = |a_{ii}-\alpha|$ (the other case can be considered analogically). Obviously, $|a_{ii}-\alpha|$ is the length of the interval $(A,B)$ between the points $A$ with the coordinates $(a_{ii}, 0)$ and $B$ with the coordinates $(\alpha, 0)$ (see Figure 1). Denote $C$ the point with the coordinates $(a_{ii}, y(a_{ii}))$. Obviously, $C \in \partial({\mathfrak D})$, and $y(a_{ii}) > 0$ is the length of the interval $(A,C)$. Now, consider the right triangle $ABC$ with the apexes $A,B,C \in \overline{{\mathfrak D}}$. Since $\mathfrak D$ is convex, its closure $\overline{{\mathfrak D}}$ is also convex and it is easy to see that all the interior points of $ABC$ belongs to $\mathfrak D$. Simple geometric observation shows that
$$r(a_{ii}) = \frac{\|AB\|\|BC\|}{\sqrt{\|AB\|^2 +\|BC\|^2}} =\frac{\|AB\|\|BC\|}{\|CA\|} = \|AH\|, $$
is the length of the altitude $AH$ to the hypotenuse $CA$ in the right triangle $ABC$. Hence we obtain the inclusion
 \begin{equation}\label{iincl}D(a_{ii}, r(a_{ii}))\cap ABC \subset \overline{{\mathfrak D}}.\end{equation}
Now consider the point $B^*$ with the coordinates $(2a_{ii} - \alpha, 0)$. The triangle $AB^*C$, that is the symmetric reflection of the triangle $ABC$ with respect to the line $AC$. Since $|a_{ii} - \beta| \geq |a_{ii} - \alpha|$, we obtain $AB^*C \subset \overline{\mathfrak D}$, and, as it follows,
\begin{equation}\label{iiincl} D(a_{ii}, r(a_{ii}))\cap AB^*C \subset \overline{{\mathfrak D}}.\end{equation}
 Inclusion \eqref{iincl} together with \eqref{iiincl} imply that the upper half of the disk $D(a_{ii}, r(a_{ii}))$ belongs to $\overline{{\mathfrak D}}$.
 Since the region $\mathfrak D$ is symmetric with respect to the real axis, we obtain the inclusion $D(a_{ii}, r(a_{ii})) \subset \overline{{\mathfrak D}}$. Thus, any closed disc, centered at the point $A(a_{ii}, 0)$ of a smaller radius, belongs to $\mathfrak D$.

\begin{figure}[h]
\center{\includegraphics[scale=0.5]{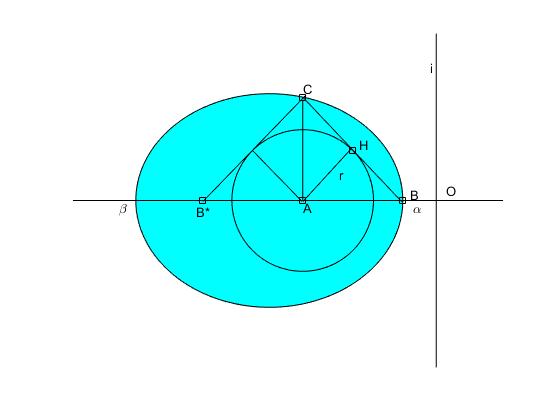}}
\caption{$A(a_{ii}, 0)$, $B(\alpha, 0)$, $C(a_{ii}, y(a_{ii}))$, $B^*(2a_{ii} - \alpha, 0)$; $AH\perp BC$.}
\end{figure}

\item[\rm Case II.] Let both $\alpha$ and $\beta$ be infinite (in this case, $\mathfrak D$ is a horizontal stripe). Then
 $$r(a_{ii}) = y(a_{ii}):=c. $$
 Here, $c$ denotes the length of the interval $(A,B)$ between the point $A$ with the coordinates $(a_{ii}, 0)$ and the point $B$ with the coordinates $(a_{ii}, c)$ (see Figure 2). Obviously, $D(a_{ii}, c) \subset \overline{{\mathfrak D}}$.
 \begin{figure}[h]
\center{\includegraphics[scale=0.5]{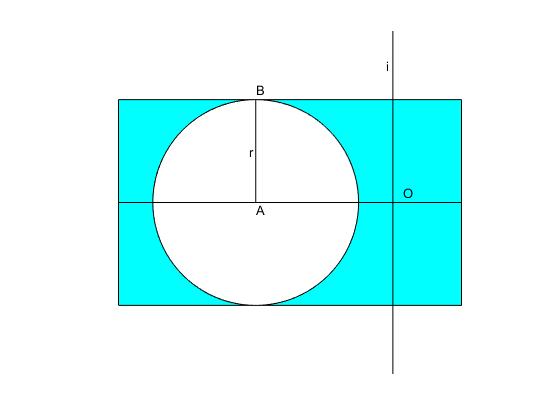}}
\caption{$A(a_{ii}, 0)$, $B(a_{ii}, c)$; $r = r(a_{ii})= c$.}
\end{figure}
\item[\rm Case III.] Let $\mathfrak D$ be unbounded from above (in this case, $\mathfrak D$ is a vertical stripe or a halfplane). Without loss the generality of the reasoning, we may assume that it is a vertical stripe, i.e. both $\alpha$ and $\beta$ are finite. Then
$$r(a_{ii}) = \min(|a_{ii}-\alpha|,|a_{ii}-\beta|). $$ Let $r(a_{ii}) = |a_{ii}-\alpha| = \alpha - a_{ii}$. Consider the points $A$ with the coordinates $(a_{ii},0)$, $B$ with coordinates $(\alpha, 0)$ and $B'$ with the coordinates $(2a_{ii} - \alpha, 0)$. Since $|a_{ii}-\alpha| = \alpha - a_{ii} \leq |a_{ii}-\beta| = a_{ii}-\beta$, we obtain $\beta \leq 2a_{ii} - \alpha$. Thus the interval $B'B$ lies entirely in $\overline{\mathfrak D}$, and it is easy to see that the circle $D(a_{ii}, \alpha - a_{ii})$ lies entirely in $\overline{{\mathfrak D}}$ (see Figure 3).
\begin{figure}[h]
\center{\includegraphics[scale=0.5]{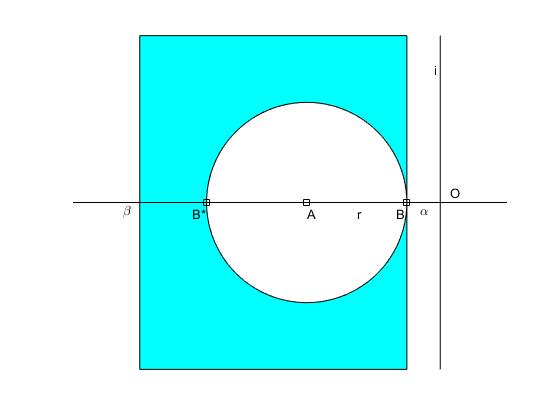}}
\caption{$A(a_{ii}, 0)$, $B(\alpha, 0)$, $B^*(2a_{ii} - \alpha, 0)$; $r = r(a_{ii}) = |a_{ii}-\alpha|$.}
\end{figure}
\end{enumerate} $\square$

\begin{corollary} Let ${\mathfrak D} \subset {\mathbb C}$ be an LMI region. If $0$ does not belong to $\mathfrak D$, then diagonally $\mathfrak D$-dominant matrices are necessarily nonsingular.
\end{corollary}
{\bf Proof.} Let $\mathbf A$ be  diagonally $\mathfrak D$-dominant. Then, by Theorem \ref{mmain}, $\sigma(\mathbf A) \subset {\mathfrak D}$. Thus $0$ does not belong to $\sigma(\mathbf A)$. $\square$

For some particular cases, we obtain preservation of $\mathfrak D$-stability under multiplication by a diagonal matrix.

\begin{theorem}\label{par} Given an unbounded LMI region $\mathfrak D$ such that ${\mathfrak D}\cap{\mathbb R} = {\mathbb R}_-$. If a matrix ${\mathbf A} \subset {\mathcal M}^{n \times n}$ is diagonally ${\mathfrak D}$-dominant then it is $\mathfrak D$-stable. Moreover, every matrix of the form ${\mathbf D}{\mathbf A}$, where ${\mathbf D} \in {\mathcal D}^+_{(0,1]}$ is also $\mathfrak D$-stable.
\end{theorem}
{\bf Proof}. Let $\mathbf A$ be diagonally ${\mathfrak D}$-dominant. Then, applying Theorem \ref{mmain}, we obtain that it is $\mathfrak D$-stable. Applying Lemma \ref{par} to the matrix ${\mathbf D}{\mathbf A}$, where ${\mathbf D}$ is an arbitrary matrix from ${\mathcal D}^+_{(0,1]}$, we obtain that ${\mathbf D}{\mathbf A}$ is also diagonally ${\mathfrak D}$-dominant, and, by Theorem \ref{mmain}, also $\mathfrak D$-stable. $\square$

\begin{theorem}\label{mainth} Let an unbounded LMI region $\mathfrak D$ possesses the characteristics $\theta_0 > 0$ and $x_{max} = \beta < +\infty$. Let a matrix ${\mathbf A} \in {\mathcal M}^{n \times n}$ be diagonally ${\mathbb C}_{\theta}^{\alpha}$-dominant with $0 < \theta \leq \theta_0$ and $\alpha \leq \min(0, x_{max})$. Then $\mathbf A$ is ${\mathfrak D}$-stable, and, moreover, $({\mathfrak D}, D)$-stable.
\end{theorem}

{\bf Proof.} First, let us prove that if $\mathbf A$ is diagonally ${\mathbb C}_{\theta}^{\alpha}$-dominant with $0 < \theta \leq \theta_0$ and $\alpha \leq \min(0, x_{max})$, then $\mathbf A$ is ${\mathfrak D}$-stable. Indeed, from the inequalities $0 < \theta \leq \theta_0$ and $\alpha \leq x_{max}$ and Inclusion \eqref{sect}, we get the inclusions
$${\mathbb C}_{\theta}^{\alpha} \subseteq {\mathbb C}_{\theta_0}^{x_{max}} \subseteq {\mathfrak D}.$$
Thus, by Lemma \ref{inclu}, if $\mathbf A$ is diagonally ${\mathbb C}_{\theta}^{\alpha}$-dominant then it is diagonally $\mathfrak D$-dominant. Applying Theorem \ref{mmain}, we obtain that $\mathbf A$ is $\mathfrak D$-stable.

Now let us prove that $\mathbf A$ is $({\mathfrak D}, D)$-stable. Here, we consider the following two cases.

\begin{enumerate}
\item[\rm Case 1.] Let $0 \in \overline{{\mathfrak D}}$. In this case, $x_{max} \geq 0$ and we get $\alpha \leq 0 \leq x_{max}$. Thus, we obtain
     $${\mathbb C}_{\theta}^{\alpha} \subseteq {\mathbb C}_{\theta}^{0},$$
     and, by Lemma \ref{inclu}, $\mathbf A$ is necessarily diagonally ${\mathbb C}_{\theta}^{0}$-dominant. Applying Lemma \ref{obv}, we obtain that ${\mathbf D}{\mathbf A}$ is also diagonally ${\mathbb C}_{\theta}^{0}$-dominant for any positive diagonal matrix $\mathbf D$. From the inclusion ${\mathbb C}_{\theta}^{0} \subseteq {\mathfrak D}$ and Lemma \ref{inclu},
     we obtain that ${\mathbf D}{\mathbf A}$ is diagonally $\mathfrak D$-dominant for every positive diagonal $\mathbf D$. Then, applying Theorem \ref{mmain}, we get that ${\mathbf D}{\mathbf A}$ is ${\mathfrak D}$-stable for every positive diagonal matrix ${\mathbf D}$. Hence, according to Definition 5, $\mathbf A$ is $({\mathfrak D}, D)$-stable.
\item[\rm Case 2.] Let $0 \in {\mathbb C} \setminus \overline{{\mathfrak D}}$. In this case, $\alpha \leq x_{max} < 0$. Applying Lemma \ref{hyp}, we obtain that if ${\mathbf A}$ is diagonally ${\mathbb C}^{\alpha}_{\theta}$-dominant, then ${\mathbf D}{\mathbf A}$ is also diagonally ${\mathbb C}^{\alpha}_{\theta}$-dominant for any ${\mathbf D} \in {\mathcal D}_{\geq 1}$.  From the inclusion ${\mathbb C}^{\alpha}_{\theta} \subseteq {\mathfrak D}$ and Lemma \ref{inclu},
     we obtain that ${\mathbf D}{\mathbf A}$ is diagonally $\mathfrak D$-dominant for any ${\mathbf D} \in {\mathcal D}_{\geq 1}$. Then, applying Theorem \ref{mmain}, we get that ${\mathbf D}{\mathbf A}$ is  ${\mathfrak D}$-stable for every ${\mathbf D} \in {\mathcal D}^+_{\geq 1}$. Hence, according to Definition 5', $\mathbf A$ is $({\mathfrak D}, D)$-stable.
\end{enumerate}
    $\square$

\section{Stability of second order systems}

 Consider the system of second-order differential equations
\begin{equation}\label{syst2}
\ddot{x} = {\mathbf A}\dot{x} + {\mathbf B}x, \qquad x \in {\mathbb R^n}, \qquad {\mathbf A}, {\mathbf B} \in {\mathcal M}^{n \times n}.
\end{equation}

Following \cite{NIS2}, we re-write System \eqref{syst2} formally as the system of first-order differential equations:

\begin{equation}\label{syst3}\begin{pmatrix}\ddot{x} \\ \dot{x} \end{pmatrix} = \begin{pmatrix}{\mathbf A} & {\mathbf B} \\ {\mathbf I} & {\mathbf O} \end{pmatrix}\begin{pmatrix}\dot{x} \\ x \end{pmatrix},\end{equation}
where $\mathbf I$ is an $n \times n$ identity matrix, $\mathbf O$ is an $n \times n$ zero matrix.

Clearly, the dynamics of System \eqref{syst2} is determined by the properties of an $2n \times 2n$ matrix of the form
\begin{equation}\label{C}\widetilde{{\mathbf A}}:= \begin{pmatrix}{\mathbf A} & {\mathbf B} \\
 {\mathbf I} & {\mathbf O} \\ \end{pmatrix}. \end{equation}
 Going into details, System \eqref{syst2} is asymptotically stable if and only if the matrix $\widetilde{\mathbf A}$ is stable, i.e. all its eigenvalues have negative real parts (see, for example, \cite{NIS}, \cite{NIS2}).

 Since the matrix $\widetilde{{\mathbf A}}$ obviously does not preserve the structural properties of $\mathbf A$ and $\mathbf B$, such as diagonal dominance, we will use the following lemma which describes the spectrum of $\widetilde{{\mathbf A}}$ (see \cite{NIS}, p. 187, Lemma 2, Part (a)).

 \begin{lemma}\label{7} Let ${\mathbf A}, \ {\mathbf B} \in {\mathcal M}^{n \times n}$ and the $2n \times 2n$ matrix $\widetilde{{\mathbf A}}$ be defined by \eqref{C}. Then each eigenvalue $\lambda$ of $\widetilde{{\mathbf A}}$ is a zero of the equation
 \begin{equation}\label{det}
 \det{\mathbf G}(\mu) = 0,
 \end{equation}
 where ${\mathbf G}(\mu) \in {\mathbb C}^{n \times n}$ is defined by the formula
 \begin{equation}\label{G}
 {\mathbf G}(\mu) = \mu^2{\mathbf I} - \mu{\mathbf A} - {\mathbf B}.
 \end{equation}
 \end{lemma}

Now let us recall the complex version of Lyapunov's matrix stability criterion (see, for example, \cite{GANT}, \cite{GANT2}).
    \begin{theorem}[Lyapunov]\label{lyap} A matrix ${\mathbf A} \in {\mathcal M}^{n \times n}({\mathbb C})$ is stable if and only if there exists a Hermitian positive definite matrix $\mathbf H$ such that the matrix
    $${\mathbf W}:={\mathbf H}{\mathbf A} + {\mathbf A}^{*}{\mathbf H}$$
    is negative definite.
    \end{theorem}

Here, we state and prove a sufficient condition for the stability of $\widetilde{{\mathbf A}}$, basing on the Lyapunov Theorem \ref{lyap}.

\begin{theorem}\label{diagcrit} Let ${\mathbf A}=\{a_{ij}\}_{i,j=1}^n$ and ${\mathbf B}=\{b_{ij}\}_{i,j = 1}^n$ be real $n \times n$ matrices, and the $2n \times 2n$ matrix $\widetilde{{\mathbf A}}$ be defined by \eqref{C}. Let $\mathbf B$ be diagonalizable, i.e there is an invertible matrix $\mathbf S$, such that ${\mathbf B} = {\mathbf S}{\mathbf \Lambda}_{\mathbf B}{\mathbf S}^{-1}$, where ${\mathbf \Lambda}_{\mathbf B} = {\rm diag}\{\nu_1, \ \ldots \nu_n\}$, and $\nu_1, \ \ldots \nu_n$ are the eigenvalues of $\mathbf B$. If $\nu_i < 0$ for all $i = 1, \ \ldots, n$, and the matrix ${\mathbf S}^{-1}{\mathbf A}{\mathbf S}$ is diagonally stable, then $\widetilde{{\mathbf A}}$ is stable. In particular, if $\mathbf B$ is negative diagonal and $\mathbf A$ is diagonally stable, then $\widetilde{{\mathbf A}}$ is stable.
\end{theorem}
{\bf Proof.} Let ${\mathbf G}(\lambda) = \lambda^2 - \lambda{\mathbf A} - {\mathbf B}$. We consider the matrix $$\widetilde{{\mathbf G}}(\lambda) = -\frac{1}{\lambda}{\mathbf S}^{-1}{\mathbf G}(\lambda){\mathbf S} = -\lambda{\mathbf I} + {\mathbf S}^{-1}{\mathbf A}{\mathbf S} + \frac{1}{\lambda}{\mathbf \Lambda}_{\mathbf B}.$$
Since $\mathbf B$ is nonsingular, ${\mathbf G}(0)$ is also nonsingular. Hence ${\mathbf G}(\lambda)$ is nonsingular if and only if $\widetilde{{\mathbf G}}(\lambda)$ is nonsingular. To establish non-singularity of $\widetilde{{\mathbf G}}(\lambda)$ for all $\lambda$ with ${\rm Re}(\lambda) \geq 0$, it is enough to prove its stability for all $\lambda$ with ${\rm Re}(\lambda) \geq 0$. From Lyapunov Theorem \ref{lyap}, we obtain, that it is necessary and sufficient for stability, to find for each $\lambda \in {\mathbb C}$ with ${\rm Re}(\lambda) \geq 0$, a Hermitian positive definite matrix ${\mathbf H}(\lambda)$ such that the matrix
    $${\mathbf W}(\lambda):={\mathbf H}(\lambda)\widetilde{{\mathbf G}}(\lambda) + (\widetilde{{\mathbf G}}(\lambda))^*{\mathbf H}(\lambda)$$
    would be negative definite.

Since ${\mathbf S}^{-1}{\mathbf A}{\mathbf S}$ is diagonally stable, it satisfies Lyapunov theorem with some positive diagonal matrix ${\mathbf D}_{\mathbf A}$, i.e the matrix
$${\mathbf W}_{\mathbf A}:={\mathbf D}_{\mathbf A}({\mathbf S}^{-1}{\mathbf A}{\mathbf S}) + ({\mathbf S}^{-1}{\mathbf A}{\mathbf S})^T{\mathbf D}_{\mathbf A},$$
    is negative definite.

Putting ${\mathbf H}(\lambda) := {\mathbf D}_{\mathbf A}$, we obtain
$${\mathbf W}(\lambda):= {\mathbf D}_{\mathbf A}\widetilde{{\mathbf G}}(\lambda) + (\widetilde{{\mathbf G}}(\lambda))^* {\mathbf D}_{\mathbf A}=$$
$${\mathbf D}_{\mathbf A}(-\lambda{\mathbf I} + {\mathbf S}^{-1}{\mathbf A}{\mathbf S} + \frac{1}{\lambda}{\mathbf \Lambda}_{\mathbf B}) + (-\lambda{\mathbf I} + {\mathbf S}^{-1}{\mathbf A}{\mathbf S} + \frac{1}{\lambda}{\mathbf \Lambda}_{\mathbf B})^*{\mathbf D}_{\mathbf A}=$$
$$-2{\rm Re}(\lambda){\mathbf D}_{\mathbf A} + {\mathbf D}_{\mathbf A}({\mathbf S}^{-1}{\mathbf A}{\mathbf S}) + ({\mathbf S}^{-1}{\mathbf A}{\mathbf S})^T{\mathbf D}_{\mathbf A} + \frac{1}{\lambda}{\mathbf D}_{\mathbf A}{\mathbf \Lambda}_{\mathbf B} + \frac{1}{\overline{\lambda}}{\mathbf \Lambda}_{\mathbf B}{\mathbf D}_{\mathbf A} =$$
$$-2{\rm Re}(\lambda){\mathbf D}_{\mathbf A} + {\mathbf W}_{\mathbf A} + \frac{2{\rm Re}(\lambda)}{|\lambda|^2}{\mathbf D}_{\mathbf A}{\mathbf \Lambda}_{\mathbf B},$$
which is negative definite whenever ${\rm Re}(\lambda) \geq 0$. $\square$

A number of classes of diagonally stable matrices is considered in the literature (see, for example, \cite{KU2}, \cite{KAB}). Let us recall, that a matrix ${\mathbf A} \in {\mathcal M}^{n \times n}$ is diagonally stable if it satisfies one of the following conditions:
\begin{enumerate}
\item[\rm 1.] Its symmetric part $\frac{{\mathbf A} +{\mathbf A}^T}{2}$ is negative definite;
\item[\rm 2.] $-\mathbf A$ is an $M$-matrix, or, equivalently, $-\mathbf A$ is a stable $Z$-matrix;
\item[\rm 3.] $n = 2$ and ${\mathbf A} = \{a_{ij}\}_{i,j=1}^2$ satisfies the conditions $a_{11}<0$, $a_{22} < 0$ and $\det({\mathbf A}) > 0$;
\item[\rm 4.] ${\mathbf A}$ is a normal stable matrix;
\item[\rm 5.] $-{\mathbf A}$ is a triangular $P$-matrix or, equivalently, $\mathbf A$ is a triangular stable matrix;
\item[\rm 6.] $\mathbf A$ is an NDD matrix;
\item[\rm 7.] $-{\mathbf A}$ is a tridiagonal $P$-matrix;
\item[\rm 8.] $-{\mathbf A}$ is a nonsingular $H_+$-matrix.
\end{enumerate}

Thus, we obtain the following corollaries (see \cite{NIS}, Theorem 2 and Theorem 1, Parts (a) and (b)).
\begin{corollary}[\cite{NIS}]\label{Crit1}
Let ${\mathbf A}=\{a_{ij}\}_{i,j=1}^n$ and ${\mathbf B}=\{b_{ij}\}_{i,j = 1}^n$ be real $n \times n$ matrices, and the $2n \times 2n$ matrix $\widetilde{{\mathbf A}}$ be defined by \eqref{C}. Let $a_{ii} < 0$ and $b_{ii} < 0$ for all $i = 1, \ \ldots, \ n$. If, in addition, $b_{ij} = 0$ for all $i \neq j$ (i.e. ${\mathbf B}$ is negative diagonal) and ${\mathbf A}$ is NDD, then $\widetilde{{\mathbf A}}$ is stable.
\end{corollary}
{\bf Proof.} For the proof, it is enough to notice, that NDD matrices are diagonally stable (see \cite{MOY}) and to apply Theorem \ref{diagcrit}. $\square$

\begin{corollary}[\cite{NIS}]\label{Crit11}
Let ${\mathbf A}=\{a_{ij}\}_{i,j=1}^n$ and ${\mathbf B}=\{b_{ij}\}_{i,j = 1}^n$ be real $n \times n$ matrices, and the $2n \times 2n$ matrix $\widetilde{{\mathbf A}}$ be defined by \eqref{C}. Let $a_{ii} < 0$ and $b_{ii} < 0$ for all $i = 1, \ \ldots, \ n$. Let, in addition, $b_{ij} = 0$ for all $i \neq j$ (i.e. ${\mathbf B}$ be negative diagonal). Then,
\begin{enumerate}
\item[\rm 1.] If $n = 1$, then $\widetilde{{\mathbf A}}$ is stable;
\item[\rm 2.] If $n = 2$ and $\mathbf A$ is stable then $\widetilde{{\mathbf A}}$ is stable.
\end{enumerate}
\end{corollary}
{\bf Proof.} For the proof, it is enough to notice, that for the case $n = 2$ stability together with negativity of the principal diagonal entries imply diagonal stability, and to apply Theorem \ref{diagcrit}. $\square$

\section{${\mathbb C}^-_{\alpha}$-stability of matrices with diagonally $\mathfrak D$-dominant submatrices}
\subsection{The boundaries for the real eigenvalues of $\widetilde{{\mathbf A}}$}
Let $\alpha$ be an arbitrary real number. Later on, we shall face the following question: "When $\alpha \in {\mathbb R}$ does not belong to $\sigma(\widetilde{{\mathbf A}})$?" or, equivalently: "When ${\mathbf G}(\alpha) = \alpha^2{\mathbf I} - \alpha{\mathbf A} - {\mathbf B}$ will be nonsingular?"

In the case when both of the matrices ${\mathbf A}=\{a_{ij}\}_{i,j=1}^n$ and ${\mathbf B}=\{b_{ij}\}_{i,j = 1}^n$ are upper (or lower) triangular, we have that ${\mathbf G}(\alpha)$ is also triangular and immediately obtain the following condition:
 $$\alpha^2 - \alpha a_{ii} - b_{ii} \neq 0 \qquad i = 1, \ \ldots, \ n, $$
 which implies
 $$b_{ii} \neq \alpha a_{ii} - \alpha^2 \qquad i = 1, \ \ldots, \ n. $$

 Now let us consider the case of arbitrary matrices $\mathbf A$ and $\mathbf B$.

\begin{lemma}\label{lem} Let ${\mathbf A}=\{a_{ij}\}_{i,j=1}^n$ and ${\mathbf B}=\{b_{ij}\}_{i,j = 1}^n$ be real $n \times n$ matrices, and the $2n \times 2n$ matrix $\widetilde{{\mathbf A}}$ be defined by \eqref{C}. Let one of the following conditions hold:
\begin{enumerate}
\item[\rm 1.] $\alpha \geq 0$; $b_{ii} < -\alpha^2$ for all $i = 1, \ \ldots, \ n$;
\item[\rm 2.]$\alpha < 0$; $b_{ii} < 3\alpha^2 -2a_{ii}\alpha$ for all $i = 1, \ \ldots, \ n$.
\end{enumerate}
 If, in addition, $b_{ij} = 0$ for all $i \neq j$ (i.e. ${\mathbf B}$ is negative diagonal) and ${\mathbf A}$ is ${\mathbb C}^-_{2\alpha}$-diagonally dominant, we obtain that $\alpha$ does not belong to $\sigma({\widetilde{\mathbf A}})$.
\end{lemma}

{\bf Proof.} For the proof, it is enough for us to show that ${\mathbf G}(\alpha)$ is diagonally dominant, thus nonsingular. We will consider the following cases.

{\bf Case I. $\alpha \geq 0$.} If $\alpha = 0$ then ${\mathbf G}(0) = - {\mathbf B}$ is positive diagonal, thus nonsingular. Now let $\alpha > 0$. Since $b_{ij} = 0$ for $i\neq j$ and \begin{equation}\label{c2} 2\alpha - a_{ii} > \sum_{i \neq j}|a_{ij}|\geq 0,\end{equation} we obtain
 $$\sum_{i \neq j}|g_{ij}(\alpha)| =\sum_{i \neq j}|\alpha a_{ij} + b_{ij}| = |\alpha|\sum_{i \neq j}|a_{ij}| < |\alpha(2\alpha - a_{ii})| = \ldots$$
  Since $\alpha > 0$; $2\alpha - a_{ii} >0$ and $b_{ii} < -\alpha^2$, we obtain
$$\ldots =\alpha(2\alpha - a_{ii}) = \alpha^2 - \alpha a_{ii} + \alpha^2 < \alpha^2 - \alpha a_{ii} - b_{ii} =|g_{ii}(\alpha)|. $$
The matrix ${\mathbf G}(\alpha)$ is diagonally dominant, thus nonsingular.

{\bf Case II. $\alpha < 0$.} In this case, $\alpha(2\alpha - a_{ii})< 0$.
Thus we obtain
$$|g_{ij}(\alpha)| = |\alpha|\sum_{i \neq j}|a_{ij}| < |\alpha(2\alpha - a_{ii})| = $$ $$\alpha a_{ii} - 2\alpha^2  = \alpha^2 - \alpha a_{ii} - 3\alpha^2 + 2\alpha a_{ii}< $$
$$ \alpha^2 - \alpha a_{ii} - b_{ii} = |g_{ii}(\alpha)|.$$
$\square$

In the sequel, we will need the following technical results.

 \begin{lemma}\label{tech} Let $\alpha, a, b \in {\mathbb R}$. Any of the conditions
 \begin{enumerate}
\item[\rm 1.] $\alpha \geq 0$; $b < -\alpha^2$;
\item[\rm 2.]$\alpha < 0$; $b < 3\alpha^2 -2a\alpha$.
\end{enumerate}
 together with the inequality $a < 2\alpha$ implies
 \begin{equation}\label{fort}
 b < \alpha(\alpha - a).
 \end{equation}
 \end{lemma}
{\bf Proof.} Condition 1 together with $a < 2\alpha$ implies:
$$2\alpha - a >0; \qquad \alpha(2\alpha - a) \geq 0; $$
$$b < -\alpha^2 \leq -\alpha^2 + \alpha(2\alpha - a) = \alpha^2 - \alpha a_. $$
Condition 2 together with $a < 2\alpha$ implies:
$$2\alpha - a >0; \qquad \alpha(2\alpha - a)< 0; $$
 $$b< 3\alpha^2 -2a\alpha =  \alpha^2 - \alpha a +\alpha(2\alpha - a) < \alpha^2 - \alpha a.$$
 $\square$

\subsection{ Triangular matrices.}
Here, we study the conditions for the localization of the spectrum of $\widetilde{{\mathbf A}}$ inside the shifted left half-plane ${\mathbf C}^-_{\alpha}$.
 First, let us consider the case, when ${\mathbf A} = \{a_{ij}\}_{i,j=1}^n$ and ${\mathbf B} = \{b_{ij}\}_{i,j = 1}^n$ both are lower (or upper) triangular. In this case, the matrix ${\mathbf G}(\lambda) = \lambda^2{\mathbf I} - \lambda{\mathbf A} - {\mathbf B}$ is also lower (respectively, upper) triangular. Obviously, $$\det{\mathbf G}(\lambda) = \prod_{i=1}^n g_{ii}(\lambda) = \prod_{i=1}^n(\lambda^2 -\lambda a_{ii} -b_{ii}).$$
Thus, the $2n$ solutions of the equation $\det{\mathbf G}(\lambda) = 0$ are exactly the solutions of $n$ quadratic equations $\lambda^2 -\lambda a_{ii} -b_{ii} = 0$, $i = 1, \ \ldots, \ n$.

Solving these quadratic equations, we obtain
$$\lambda_{2i-1, 2i}= \frac{a_{ii}\pm \sqrt{a^2_{ii} + 4b_{ii}}}{2} \qquad i = 1, \ \ldots, \ n.$$

Now let us analyze the conditions when all these solutions belong to the shifted half-plane ${\mathbb C}^-_{\alpha}$, $\alpha \in {\mathbb R}$. Here, we have the following two cases.

\begin{enumerate}
\item[\rm 1.] $D = a^2_{ii} + 4b_{ii} \geq 0$. In this case, we have two (probably, coinciding) real solutions $\lambda_{2i-1, 2i}$, that should satisfy the inequalities $\lambda_{2i-1, 2i} < \alpha$. Hence we obtain the following conditions:
    $$a_{ii}\pm \sqrt{a^2_{ii} + 4b_{ii}} < 2\alpha;$$
    $$\pm \sqrt{a^2_{ii} + 4b_{ii}} < 2\alpha - a_{ii} \ \Rightarrow \ 2\alpha - a_{ii} > 0;$$
    $$a^2_{ii} + 4b_{ii} < 4\alpha^2 + a^2_{ii} - 4\alpha a_{ii}; $$
    $$b_{ii} < \alpha^2 - \alpha a_{ii}. $$
Summarizing the above conditions, we obtain the following set of inequalities:
$$\frac{a_{ii}}{2} < \alpha; \qquad -\frac{a_{ii}^2}{4} \leq b_{ii} < \alpha^2 - \alpha a_{ii}. $$
\item[\rm 2.] $D = a^2_{ii} + 4b_{ii} < 0$. In this case, we have two non-real solutions $\lambda_{2i-1, 2i}$, given by the formula
$$\lambda_{2i-1, 2i}= \frac{a_{ii}\pm i\sqrt{-(a^2_{ii} + 4b_{ii})}}{2}.$$
In this case, ${\rm Re}(\lambda_{2i-1, 2i}) = \frac{a_{ii}}{2}$, and we immediately obtain the conditions $$\frac{a_{ii}}{2} < \alpha; \qquad  b_{ii} < -\frac{a_{ii}^2}{4}. $$
\end{enumerate}

Summarizing both of the cases, we obtain the following statement.
\begin{theorem}
Let ${\mathbf A}=\{a_{ij}\}_{i,j=1}^n$ and ${\mathbf B}=\{b_{ij}\}_{i,j = 1}^n$ be both real $n \times n$ upper (or lower) triangular matrices, and the $2n \times 2n$ matrix $\widetilde{{\mathbf A}}$ be defined by \eqref{C}. Then $\widetilde{{\mathbf A}}$ is ${\mathbb C}^-_{\alpha}$-stable if and only if $\frac{a_{ii}}{2} < \alpha$ and $b_{ii} < \alpha^2 - \alpha a_{ii}$.
\end{theorem}
{\bf Proof.} The proof follows from the above reasoning. $\square$

\subsection{Case of arbitrary matrices}
Here, we consider the following sufficient condition, based on the property of diagonal ${\mathbb C}^-_{\alpha}$-dominance.

\begin{theorem}\label{T2}
Let ${\mathbf A}=\{a_{ij}\}_{i,j=1}^n$ and ${\mathbf B}=\{b_{ij}\}_{i,j = 1}^n$ be real $n \times n$ matrices, and the $2n \times 2n$ matrix $\widetilde{{\mathbf A}}$ be defined by \eqref{C}. Let $b_{ii} < -\alpha^2$ if $\alpha \geq 0$ and $b_{ii} < 3\alpha^2 -2a_{ii}\alpha$ if $\alpha < 0$ for all $i = 1, \ \ldots, \ n$. If, in addition, $b_{ij} = 0$ for all $i \neq j$ (i.e. ${\mathbf B}$ is negative diagonal) and ${\mathbf A}$ is diagonally ${\mathbb C}^-_{2\alpha}$-dominant, then $\widetilde{{\mathbf A}}$ is ${\mathbb C}^-_{\alpha}$-stable.
\end{theorem}
{\bf Proof.} Applying Lemma \ref{lem}, we exclude the case $\lambda = \alpha$, and examine the matrix $\widetilde{{\mathbf G}}(\lambda) = - \frac{1}{\lambda - \alpha}{\mathbf G}(\lambda)$. For its principal diagonal entries, we have the equality
 $$\widetilde{g}_{ii}(\lambda) = \frac{ -\lambda^2 + \lambda a_{ii} + b_{ii}}{\lambda - \alpha} = \ldots$$
 Applying the Taylor expansion formula, we obtain
  $$\ldots = \frac{ (-\alpha^2 + \alpha a_{ii} + b_{ii}) + (a_{ii} - 2\alpha)(\lambda - \alpha) - (\lambda - \alpha)^2}{\lambda - \alpha} = $$ $$  \frac{ -\alpha^2 + \alpha a_{ii} + b_{ii}}{\lambda - \alpha} +  (a_{ii} - 2\alpha) - (\lambda - \alpha).$$
  Let $\lambda = x+iy$, then $\lambda -\alpha= (x-\alpha)+iy$, $\frac{1}{\lambda-\alpha} = \frac{x-\alpha}{(x-\alpha)^2 + y^2} -i\frac{y}{(x-\alpha)^2 + y^2}$. Then the following estimates hold:
 $$|\widetilde{g}_{ii}(\lambda)| = |\frac{ -\alpha^2 + \alpha a_{ii} + b_{ii}}{\lambda - \alpha} +  (a_{ii} - 2\alpha) - (\lambda - \alpha)| \geq $$ $$ |{\rm Re}(\frac{ -\alpha^2 + \alpha a_{ii} + b_{ii}}{\lambda - \alpha} +  (a_{ii} - 2\alpha) - (\lambda - \alpha))| = $$ $$ |(a_{ii} - 2\alpha) -  (x-\alpha) + \frac{(x-\alpha)(-\alpha^2 + \alpha a_{ii} + b_{ii})}{(x-\alpha)^2+y^2}| = \ldots $$
 From diagonal ${\mathbb C}^-_{2\alpha}$-dominance (see Definition 6) of $\mathbf A$ we get the estimate $ a_{ii} < 2\alpha$. Applying Lemma \ref{tech}, we obtain the estimate $-\alpha^2 + \alpha a_{ii} + b_{ii} < 0$. Then for $x \geq \alpha$, all the summands $a_{ii} - 2\alpha$, $-(x-\alpha)$ and $\frac{(x-\alpha)(-\alpha^2 + \alpha a_{ii} + b_{ii})}{(x-\alpha)^2+y^2}$ are nonpositive. We have:
 $$\ldots = |a_{ii} - 2\alpha| + |x-\alpha| + |\frac{(x-\alpha)(-\alpha^2 + \alpha a_{ii} + b_{ii})}{(x-\alpha)^2+y^2}| \geq |a_{ii}-2\alpha|.$$

 Imposing the condition $b_{ij} = 0$ for $i \neq j$ (i.e. that the matrix $\mathbf B$ is diagonal), we obtain that the off-diagonal entries of $\widetilde{{\mathbf G}}(\lambda)$ coincide with the off-diagonal entries of $\mathbf A$. Then the condition
 $$|a_{ii} - 2\alpha| > \sum_{j\neq i}|a_{ij}| \qquad i = 1, \ \ldots, \ n.$$
  implies the strict diagonal dominance of $\widetilde{{\mathbf G}}(\lambda)$ (see Definition 1), and, by Gershgorin theorem, its invertibility. $\square$
  
  The above theorem also implies Corollary \ref{Crit1}, providing a direct proof for it by the methods similar to that used in \cite{NIS}, the proof of Theorem 2.

\subsection{The case of simultaneously triagonalizable matrices}

\begin{theorem}\label{RH} Let ${\mathbf A}=\{a_{ij}\}_{i,j=1}^n$ and ${\mathbf B}=\{b_{ij}\}_{i,j = 1}^n$ be real $n \times n$ matrices, and the $2n \times 2n$ matrix $\widetilde{{\mathbf A}}$ be defined by \eqref{C}. Let $\mathbf A$ and $\mathbf B$ can be simultaneously reduced (by similarity transformation) to their triangular forms. Denote $\sigma(\widetilde{\mathbf A}) = \{\lambda_k\}_{k = 1}^{2n}$, $\sigma(\mathbf A) = \{\mu_i\}_{i = 1}^n$, $\mu_i = a_i + ic_i$, $i = 1, \ \ldots, \ n$ and $\sigma(\mathbf B) = \{\nu_i\}_{i = 1}^n$, $\nu_i = b_i + id_i$, $i = 1, \ \ldots, \ n$. Then $\widetilde{\mathbf A}$ is ${\mathbb C}^-_{\alpha}$-stable if and only if ${\mathbf A}$ is ${\mathbb C}^-_{2\alpha}$-stable and for each $i$, $i = 1, \ \ldots, \ n$, the values $a_i, \ b_i, \ c_i, d_i$ and $\alpha$ satisfy the following inequality:
\begin{equation}\label{eqal}
(2\alpha - a_i)^2(\alpha^2 - a_i\alpha -b_i) + (2\alpha - a_i)c_i(\alpha c_i+d_i) - (\alpha c_i+d_i)^2 > 0.
\end{equation}
In particular, $\widetilde{\mathbf A}$ is stable if and only if ${\mathbf A}$ is stable and for each $i$, $i = 1, \ \ldots, \ n$, the values $a_i, \ b_i, \ c_i, d_i$ satisfy the following inequality:
\begin{equation}\label{eqal2}
a_i^2b_i + a_ic_id_i + d_i^2 < 0.
\end{equation}
\end{theorem}
{\bf Proof.} By Lemma \ref{7}, $\lambda$ is an eigenvalue of $\widetilde{{\mathbf A}}$ if and only if $\det{\mathbf G}(\lambda) = 0$, where
 $${\mathbf G}(\lambda) = \lambda^2{\mathbf I} - \lambda{\mathbf A} - {\mathbf B}.$$
 Since $\mathbf A$ and $\mathbf B$ are simultaneously triagonalizable, we get $\det{\mathbf G}(\lambda) = 0$ implies $\det({\mathbf S}{\mathbf G}(\lambda){\mathbf S}^{-1}) = 0$, where ${\mathbf S}{\mathbf A}{\mathbf S}^{-1} = {\mathbf \Lambda}_{\mathbf A}$, ${\mathbf S}{\mathbf B}{\mathbf S}^{-1} = {\mathbf \Lambda}_{\mathbf B}$, ${\mathbf \Lambda}_{\mathbf A}$ and ${\mathbf \Lambda}_{\mathbf B}$ are the triangular forms of the matrices $\mathbf A$ and $\mathbf B$, respectively. Then, denoting $\mu_i$ and $\nu_i$, $i = 1, \ \ldots, n$ the eigenvalues of the matrices $\mathbf A$ and $\mathbf B$, respectively, we obtain the equalities:
$$ 0 = \det({\mathbf S}{\mathbf G}(\lambda){\mathbf S}^{-1}) = \det(\lambda^2{\mathbf I} - \lambda{\mathbf \Lambda}_{\mathbf A} - {\mathbf \Lambda}_{\mathbf B}) = \prod_{i = 1}^n(\lambda^2 -\lambda\mu_i - \nu_i). $$
Thus the spectrum of $\widetilde{{\mathbf A}}$ coincides with the solutions of $n$ quadratic equations
 \begin{equation}\label{Eq2}\lambda^2 -\lambda\mu_i - \nu_i = 0, \qquad i = 1, \ \ldots, \ n. \end{equation}
 Consider the polynomial $$f_i(\lambda) = \lambda^2 -\lambda\mu_i - \nu_i = \lambda^2 + \lambda(-a-ic) + (-b-id),$$
 with complex coefficients $-\mu_i = -a_i-ic_i$ and $-\nu_i = -b_i-id_i$. It is easy to see, that the polynomial $f_i(\lambda)$ has all its zeroes in ${\mathbb C}^-_{\alpha}$ if and only if the polynomial $f_i(\lambda+\alpha)$ is stable. Consider
 $$\hat{f}_i(\lambda):= f_i(\lambda+\alpha) = (\lambda+\alpha)^2 -(\lambda+\alpha)\mu_i - \nu_i=$$
 $$=\lambda^2 + \lambda(2\alpha -\mu_i) +(\alpha^2-\mu_i\alpha-\nu_i)=$$
 $$\lambda^2 +\lambda((2\alpha - a_i)+i(-c_i)) + (\alpha^2 - a_i\alpha -b_i) +i(-\alpha c_i-d_i) = \lambda^2 + \lambda(\hat{a}_i+i\hat{c}_i) + (\hat{b}_i+i\hat{d}_i). $$
 Applying the generalized Routh--Hurwitz criterion (see \cite{XI}, the generalized Routh array, also see \cite{GANT}) to the polynomial $\hat{f}_i(\lambda)$, we obtain that  $\hat{f}_i(\lambda)$ is stable if and only if the following inequalities hold:
 \begin{equation}\label{D1} \Delta_1 =\hat{a} = 2\alpha - a_i > 0; \end{equation}
 \begin{equation}\label{D2} \Delta_2 = \hat{a}^2\hat{b} + \hat{a}\hat{c}\hat{d} - \hat{d}^2 > 0.\end{equation}
 Condition \eqref{D1} exactly means $a_i = {\rm Re}(\mu_i)< 2\alpha$ for each eigenvalue $\mu_i$ of $\mathbf A$. Substituting the expressions for the coefficients of $f_i(\lambda)$ to Condition \eqref{D2}, we obtain
 $$\Delta_2 = (2\alpha - a_i)^2(\alpha^2 - a_i\alpha -b_i) + (2\alpha - a_i)c_i(\alpha c_i+d_i) - (\alpha c_i+d_i)^2 > 0.$$
 Substituting $\alpha=0$ to Conditions \eqref{D1} and \eqref{D2}, we obtain stability of $\mathbf A$ and Inequality \eqref{eqal2}. $\square$

\begin{corollary}[\cite{NIS}, Theorem 3]
Let ${\mathbf A}=\{a_{ij}\}_{i,j=1}^n$ be a real $n \times n$ matrix with  $a_{ii} < 0$ for all $i = 1, \ \ldots, \ n$ and let the $2n \times 2n$ matrix $\widetilde{{\mathbf A}}$ be defined as follows: $$\widetilde{{\mathbf A}}= \begin{pmatrix}{\mathbf A} & b{\mathbf I} \\
 {\mathbf I} & {\mathbf O} \\ \end{pmatrix}, $$ where $ b < 0$. Then $\widetilde{{\mathbf A}}$ is stable if and only if $\mathbf A$ is stable.
\end{corollary}
{\bf Proof.} For the proof, it is enough to put $d_i = 0$ to Inequality \eqref{eqal2}.

The following corollaries give stronger versions of \cite{NIS}, Theorem 4 and Corollary 1.
\begin{corollary}\label{2type}
Let ${\mathbf B}=\{b_{ij}\}_{i,j = 1}^n$ be a real $n \times n$ matrix, and let the $2n \times 2n$ matrix $\widetilde{{\mathbf A}}$ be defined by $$\widetilde{{\mathbf A}}= \begin{pmatrix} a{\mathbf I} & {\mathbf B} \\
 {\mathbf I} & {\mathbf O} \\ \end{pmatrix},$$ where $a < 0$.
Then $\widetilde A$ is stable if and only if $\mathbf B$ is stable, and, in addition, $|a| > \frac{|d_i|}{\sqrt{-b_i }}$ for each eigenvalue $\nu_i = b_i +id_i$ of $\mathbf B$.
\end{corollary}
{\bf Proof.} Putting $c_i = 0$ to Inequality \eqref{eqal2}, we obtain the inequality $a^2b_i + d_i^2< 0$, which obviously implies $b_i < 0$. Taking into account $b_i < 0$, we obtain $a^2 > \frac{d_i^2}{-b_i}$ which leads to $|a| > \frac{|d_i|}{\sqrt{-b_i }}$.

\begin{corollary}
Let ${\mathbf B}=\{b_{ij}\}_{i,j = 1}^n$ be a real $n \times n$ matrix, and let the $2n \times 2n$ matrix $\widetilde{{\mathbf A}}$ be defined by $$\widetilde{{\mathbf A}}= \begin{pmatrix} a{\mathbf I} & {\mathbf B} \\
 {\mathbf I} & {\mathbf O} \\ \end{pmatrix},$$ where $a < 0$. Let all the eigenvalues of ${\mathbf B}$ be real, then $\widetilde{{\mathbf A}}$ is stable if and only if all the eigenvalues of $\mathbf B$ are negative.
\end{corollary}
{\bf Proof.} For the proof, it is enough to put $c_i = 0$ and $d_i = 0$ to Inequality \eqref{eqal2}.

Consider the following particular case of Theorem \ref{RH} for $\alpha \in {\mathbb R}$.

\begin{corollary}
Let ${\mathbf B}=\{b_{ij}\}_{i,j = 1}^n$ be a real $n \times n$ matrix, and let the $2n \times 2n$ matrix $\widetilde{{\mathbf A}}$ be defined by $$\widetilde{{\mathbf A}}= \begin{pmatrix} a{\mathbf I} & {\mathbf B} \\
 {\mathbf I} & {\mathbf O} \\ \end{pmatrix}.$$ where $a < 2\alpha$. Let all the eigenvalues of ${\mathbf B}$ be real, then $\widetilde{{\mathbf A}}$ is ${\mathbb C}^-_{\alpha}$-stable if and only if all the eigenvalues $b_i$, $i = 1, \ \ldots, \ n$ of $\mathbf B$ satisfy $b_i < \alpha^2 - a\alpha$.
\end{corollary}
{\bf Proof.} For the proof, it is enough to put $c_i = 0$ and $d_i = 0$ to Inequality \eqref{eqal}.

\section{Stability and transient response properties of perturbed second-order systems}

\subsection{Perturbations of second-order systems} Let us provide the following definitions, basing on the definitions introduced in \cite{KOK}, for linear mechanical systems.

{\bf Definition 10}. System \eqref{syst2} is called {\it (multiplicative) $D$-stable} if the system
\begin{equation}\label{systper1}
\ddot{x} = {\mathbf D}({\mathbf A}\dot{x} + {\mathbf B}x)
\end{equation}
is asymptotically (Lyapunov) stable for every positive diagonal matrix $\mathbf D$.

The localization of $\sigma(\widetilde{{\mathbf A}})$ in the shifted half-plane ${\mathbb C}^-_{-\alpha}$, $\alpha > 0$, guarantees the minimal decay rate $\alpha$ (see, e.g. \cite{GUJU}). Basing on this fact, we introduce the following definition.

{\bf Definition 11}. System \eqref{syst2} is called {\it relatively $D$-stable with the minimal decay rate $\alpha > 0$} if the minimal decay rate of the perturbed system
$$\ddot{x} = {\mathbf D}({\mathbf A}\dot{x} + {\mathbf B}x)$$
is bigger than $\alpha$ for every diagonal matrix ${\mathbf D} \in {\mathcal D}^+_{\geq 1}$.

In the sequel, we shall also consider the perturbations of System \eqref{syst2} of Form I:
\begin{equation}\label{Form1}
\ddot{x} = {\mathbf D}{\mathbf A}\dot{x} + {\mathbf B}x, \qquad x \in {\mathbb R^n},
\end{equation}
and of Form II:
\begin{equation}\label{Form2}
\ddot{x} = {\mathbf A}\dot{x} + {\mathbf D}{\mathbf B}x, \qquad x \in {\mathbb R^n},
\end{equation}
where ${\mathbf D}$ is $n \times n$ positive diagonal matrix (or an arbitrary matrix from a given subclass of ${\mathcal D}^+$).

\subsection{Stability of perturbed second-order systems}
Let us provide the following condition sufficient for $D$-stability of System \eqref{syst2}, basing on Theorem \ref{diagcrit}.

\begin{theorem}\label{perst}
 Let a system of second-order ODE be defined by \eqref{syst2}. If $\mathbf A$ is diagonally stable and $\mathbf B$ is negative diagonal, then System \eqref{syst2} is $D$-stable. \end{theorem}

{\bf Proof.}
 Clearly, $D$-stability of System \eqref{syst2} is determined by the stability of a $2n \times 2n$ perturbed matrix of the form
 $$\widetilde{{\mathbf A}}_{\mathbf D} = \begin{pmatrix} {\mathbf D}{\mathbf A} & {\mathbf D}{\mathbf B} \\
 {\mathbf I} & {\mathbf O} \\ \end{pmatrix},$$
 where $\mathbf D$ is an arbitrary $n \times n$ positive diagonal matrix. Diagonal stability of $\mathbf A$ implies that ${\mathbf D}{\mathbf A}$ is also diagonally stable for every positive diagonal matrix $\mathbf D$. Obviously, ${\mathbf D}{\mathbf B}$ is negative diagonal whenever $\mathbf B$ is negative diagonal. Thus applying Theorem \ref{diagcrit}, we obtain the stability of $\widetilde{{\mathbf A}}_{\mathbf D}$ for every positive diagonal $\mathbf D$.
 $\square$

Note, that in the statement of Theorem \ref{perst}, the matrix $\mathbf A$ is not assumed to be symmetric.

Basing on Corollary \ref{2type}, we prove a stronger version of the result from \cite{KU2} (see \cite{KU2}, p. 654, Theorem 3).

\begin{theorem} Let a system of second-order ODE be defined by \eqref{syst2}. Let ${\mathbf A} = a{\mathbf I}$ with $a < 0$. Then System \eqref{syst2} preserves stability under perturbations of Form II (see \eqref{Form2}) if and only if $\mathbf B$ is $(P(a), D)$-stable, where $P(a)$ is the parabolic stability region, defined as follows:
$$P(a) = \{z=x+iy \in {\mathbb C}: y^2 < -a^2x\}. $$
\end{theorem}

{\bf Proof}. For the proof, it is enough to apply Corollary \ref{2type} to the perturbed matrix $$\widetilde{{\mathbf A}}^{II}_{\mathbf D} = \begin{pmatrix} a{\mathbf I} & {\mathbf D}{\mathbf B} \\
 {\mathbf I} & {\mathbf O} \\ \end{pmatrix},$$
 where $\mathbf D$ is an $n \times n$ positive diagonal matrix. $\square$

Note, that the class of $(P(a), D)$-stable matrices contains $D$-negative matrices (a matrix ${\mathbf A} \in {\mathcal M}^{n \times n}$ is called $D$-negative, if its spectrum $\sigma({\mathbf A})$ is real and negative and $\sigma({\mathbf D}{\mathbf A})$ is also real and negative for any positive diagonal matrix $\mathbf D$), that are studied in \cite{BAO}, see also \cite{KU2}.

Theorem \ref{par} provides the following sufficient condition for the stability of perturbed System \eqref{syst2}.

\begin{theorem} Let a system of second-order ODE be defined by \eqref{syst2}. Let ${\mathbf A} = a{\mathbf I}$ with $a < 0$ and $\mathbf B$ be diagonally $P(a)$-dominant (see Definition 9). Then System \eqref{syst2} preserves stability under perturbations of Form II (see \eqref{Form2}), where ${\mathbf D}$ is an arbitrary matrix from ${\mathcal D}^+_{(0,1]}$.
\end{theorem}

{\bf Proof.} Let us consider the perturbed matrix $$\widetilde{{\mathbf A}}^{II}_{\mathbf D} = \begin{pmatrix} a{\mathbf I} & {\mathbf D}{\mathbf B} \\
 {\mathbf I} & {\mathbf O} \\ \end{pmatrix},$$
 where ${\mathbf D} \in {\mathcal D}^+_{(0,1]}$. The matrix $\mathbf B$ is diagonally $P(a)$-dominant, hence by Theorem \ref{par}, $\mathbf B$ is $P(a)$-stable and, moreover, ${\mathbf D}{\mathbf B}$ is also $P(a)$-stable for every ${\mathbf D} \in {\mathcal D}^+_{(0,1]}$. Applying Corollary \ref{2type} to the perturbed matrix $\widetilde{{\mathbf A}}^{II}_{\mathbf D}$, we complete the proof. $\square$

\subsection{Minimal decay rate of perturbed second-order systems}

\begin{theorem}
Let a system of second-order ODE be defined by \eqref{syst2}. If $\mathbf A$ is diagonally ${\mathbb{C}^{-}_{-2\alpha}}$-dominant, $\mathbf B$ is negative diagonal and, in addition, their principal diagonal entries $a_{ii}$ and $b_{ii}$ are connected by the inequalities $b_{ii} < 2a_{ii}\alpha$, $\alpha > 0$, then System \eqref{syst2} is relatively $D$-stable with the minimal decay rate $\alpha$.
\end{theorem}

{\bf Proof.} By definition, relative $D$-stability of System \eqref{syst2} with the minimal decay rate $\alpha$ is equivalent to ${\mathbb C}^-_{-\alpha}$-stability of a $2n \times 2n$ perturbed matrix of the form
 $$\widetilde{{\mathbf A}}_{\mathbf D} = \begin{pmatrix} {\mathbf D}{\mathbf A} & {\mathbf D}{\mathbf B} \\
 {\mathbf I} & {\mathbf O} \\ \end{pmatrix},$$
 where ${\mathbf D} \in {\mathcal D}^+_{\geq 1}$. If $\mathbf A$ is diagonally ${\mathbb{C}^{-}_{-2\alpha}}$-dominant, then, by Lemma \ref{hyp}, ${\mathbf D}{\mathbf A}$ is also diagonally ${\mathbb{C}^{-}_{-2\alpha}}$-dominant for every diagonal matrix ${\mathbf D} \in {\mathcal D}^+_{\geq 1}$. Obviously, ${\mathbf D}{\mathbf B}$ is negative diagonal whenever $\mathbf B$ is negative diagonal. Finally, we need to verify, that the principal diagonal entries $d_{ii}a_{ii}$ and $d_{ii}b_{ii}$ satisfy the inequality $d_{ii}b_{ii} < 3\alpha^2 +2d_{ii}a_{ii}\alpha$, for all $i = 1, \ \ldots, n$, $d_{ii} \geq 1$. Indeed, the inequality $b_{ii} < 2a_{ii}\alpha$ implies $d_{ii}b_{ii} < 2d_{ii}a_{ii}\alpha < 3\alpha^2 +2d_{ii}a_{ii}\alpha$.

  Thus applying Theorem \ref{T2}, we obtain ${\mathbb C}^-_{-\alpha}$-stability of $\widetilde{{\mathbf A}}_{\mathbf D}$ for every positive diagonal ${\mathbf D} \in {\mathcal D}^+_{\geq 1}$.
 $\square$

 Now let us consider the cases, when the minimal decay rate of the system does not decrease under perturbations of Form I or Form II.

 \begin{theorem} Let a system of second-order ODE be defined by \eqref{syst2}. Let ${\mathbf A} = a{\mathbf I}$ with $a < 0$. Then the minimal decay rate $\alpha$ of System \eqref{syst2} does not decrease under perturbations of Form II (see \eqref{Form2}) if and only if $\mathbf B$ is $({\mathfrak D}, D)$-stable, where ${\mathfrak D}$ is the parabolic stability region, defined as follows:
$${\mathfrak D} = \{z=x+iy \in {\mathbb C}: y^2 < -(2\alpha - a)^2x + (2\alpha - a)^2(\alpha^2 -a\alpha)\}. $$
\end{theorem}

{\bf Proof}. For the proof, it is enough to apply Theorem \ref{RH} to the perturbed matrix $$\widetilde{{\mathbf A}}^{II}_{\mathbf D} = \begin{pmatrix} a{\mathbf I} & {\mathbf D}{\mathbf B} \\
 {\mathbf I} & {\mathbf O} \\ \end{pmatrix},$$
 where $\mathbf D$ is an $n \times n$ positive diagonal matrix. $\square$

The analogous study of perturbations of Type I leads to $({\mathfrak D}, D)$-stability with respect to a specific region $\mathfrak D$, bounded by some third-order curves.

\section{Stability of perturbed fractional order systems}
Consider a linear system in the following form:
\begin{equation}\label{sysf}d^{\gamma}x(t) = {\mathbf A}x(t), \end{equation}
with $0 < \gamma \leq 2$, $x(0) = x_0$. It is known (see \cite{MAT}, \cite{MAT1}, \cite{MOS}, \cite{SMF}) that System \eqref{sysf} is asymptotically stable if and only if all eigenvalues $\lambda$ of ${\mathbf A}$ satisfy the inequality
\begin{equation}\label{confrac}|\arg(\lambda)|> \gamma\dfrac{\pi}{2}.\end{equation}

\subsection{Case $0 <\gamma \leq 1$.} For $0 \leq\gamma < 1$, Condition \eqref{confrac} corresponds to ${\mathfrak D}$-stability of the system matrix $\mathbf A$ with respect to the non-convex stability region ${\mathfrak D} = {\mathbb C}\setminus\overline{{\mathbb C}^+_{\theta}}$, where $\theta = \frac{\pi\gamma}{2}$. Note, that it is not an LMI region.

 Let the dimension $n$ be even. Using the techniques, developed in the previous sections, we obtain the following conditions of matrix eigenvalue localization outside the region $\overline{{\mathbb C}^+_{\theta}}$.

 Given an $n \times n$ matrix $\mathbf A$, where $n = 2k$, we represent it in the following block form
\begin{equation}\label{ym}
{\mathbf A} = \begin{pmatrix} {\mathbf A}_{11} & {\mathbf A}_{12} \\
{\mathbf A}_{21} & {\mathbf A}_{22} \end{pmatrix},
\end{equation}
where ${\mathbf A}_{11}$, ${\mathbf A}_{12}$, ${\mathbf A}_{21}$, ${\mathbf A}_{22}$ are $k \times k$ matrices.

\begin{lemma}\label{1} Let a $2k \times 2k$ matrix ${\mathbf A}$ be defined by Equation \eqref{ym}, with $\det({\mathbf A}_{12}) \neq 0$. Then $\lambda$ is an eigenvalue of ${\mathbf A}$ if and only if $\det({\mathbf G}(\lambda)) = 0$, where
              $${\mathbf G}(\lambda) = \lambda^2{\mathbf I} - \widehat{{\mathbf A}}\lambda - \widehat{{\mathbf B}},$$
              \begin{equation}\label{A}\widehat{{\mathbf A}} = {\mathbf A}_{12}^{-1}{\mathbf A}_{11}{\mathbf A}_{12} + {\mathbf A}_{22};\end{equation}
              \begin{equation}\label{B}\widehat{{\mathbf B}} = {\mathbf A}_{21}{\mathbf A}_{12} - {\mathbf A}_{22}{\mathbf A}_{12}^{-1}{\mathbf A}_{11}{\mathbf A}_{12}.\end{equation}
\end{lemma}
{\bf Proof.} First, consider the matrix $\mathbf P$, defined as follows:
$${\mathbf P} = \begin{pmatrix}{\mathbf A}_{12}^{-1} & {\mathbf O} \\ {\mathbf O} & {\mathbf I}  \end{pmatrix}\begin{pmatrix} {\mathbf A}_{11} & {\mathbf A}_{12} \\
{\mathbf A}_{21} & {\mathbf A}_{22} \end{pmatrix}\begin{pmatrix}{\mathbf A}_{12} & {\mathbf O} \\ {\mathbf O} & {\mathbf I}  \end{pmatrix} = \begin{pmatrix} {\mathbf A}_{12}^{-1}{\mathbf A}_{11}{\mathbf A}_{12} & {\mathbf I} \\
{\mathbf A}_{21}{\mathbf A}_{12} & {\mathbf A}_{22} \end{pmatrix}.$$
Obviously, $\sigma({\mathbf P}) = \sigma({\mathbf A})$. Then, applying \cite{BER}, p. 135, Fact 2.14.13, we obtain
$$\det({\mathbf P} - \lambda{\mathbf I}) = \det\begin{pmatrix} {\mathbf A}_{12}^{-1}{\mathbf A}_{11}{\mathbf A}_{12} - \lambda{\mathbf I} & {\mathbf I} \\
{\mathbf A}_{21}{\mathbf A}_{12} & {\mathbf A}_{22} - \lambda{\mathbf I} \end{pmatrix} = $$ $$ \det(({\mathbf A}_{22} - \lambda{\mathbf I})({\mathbf A}_{12}^{-1}{\mathbf A}_{11}{\mathbf A}_{12} - \lambda{\mathbf I}) - {\mathbf A}_{21}{\mathbf A}_{12}) = $$ $$  \det(\lambda^2{\mathbf I}- \lambda({\mathbf A}_{12}^{-1}{\mathbf A}_{11}{\mathbf A}_{12} + {\mathbf A}_{22}) +{\mathbf A}_{22}{\mathbf A}_{12}^{-1}{\mathbf A}_{11}{\mathbf A}_{12} - {\mathbf A}_{21}{\mathbf A}_{12}) = $$ $$\det(\lambda^2 - \widehat{{\mathbf A}}\lambda - \widehat{{\mathbf B}}). $$
 $\square$

Note, that for the case of nonsingular $\mathbf B$ and $n:=k$, Lemma \ref{1} gives a generalization of Lemma \ref{7} and uses the same methods of the proof. Indeed, if we put $$\widetilde{{\mathbf A}} = \begin{pmatrix} {\mathbf A}_{11} & {\mathbf A}_{12} \\
{\mathbf A}_{21} & {\mathbf A}_{22} \end{pmatrix}, $$
where ${\mathbf A}_{11} = {\mathbf A}$, ${\mathbf A}_{12} = {\mathbf B}$, ${\mathbf A}_{21} = {\mathbf I}$ and ${\mathbf A}_{22} = {\mathbf O}$, Lemma \ref{1} gives us $\widehat{{\mathbf A}} = {\mathbf B}^{-1}{\mathbf A}{\mathbf B}$ and $\widehat{{\mathbf B}} = {\mathbf B}$. Hence the matrix ${\mathbf G}(\lambda)$, given by Lemma \ref{1}, is as follows:
$${\mathbf G}(\lambda) = \lambda^2{\mathbf I} - \lambda{\mathbf B}^{-1}{\mathbf A}{\mathbf B} - {\mathbf B} =$$
$$= {\mathbf B}^{-1}(\lambda^2{\mathbf I} - \lambda{\mathbf A} -  {\mathbf B}) {\mathbf B} = {\mathbf B}^{-1}\widehat{{\mathbf G}}(\lambda){\mathbf B},$$
where $\widehat{{\mathbf G}}(\lambda) = \lambda^2{\mathbf I} - \lambda{\mathbf A} -  {\mathbf B}$ is the matrix, given in the statement of Lemma \ref{7}. These two matrices are similar, thus their spectra coincide.

Now let us prove a sufficient condition for $\mathbf A$ does not have any eigenvalues in the conic sector $\overline{{\mathbb C}^+_{\theta}}$, $\theta \in (0, \frac{\pi}{4}]$.

\begin{theorem}\label{14}
Let $n = 2k$, $\theta \in (0, \frac{\pi}{4}]$ and an $n \times n$ matrix ${\mathbf A}$ be defined by \eqref{ym}, with $\det({\mathbf A}_{12}) \neq 0$. Let $\widehat{{\mathbf A}} = \{\widehat{a}_{ij}\}_{i,j = 1}^n$ and $\widehat{\mathbf B}= \{\widehat{b}_{ij}\}_{i,j = 1}^n$ be defined by \eqref{A} and \eqref{B}, respectively. If $\widehat{{\mathbf A}}$ is diagonally ${\mathbb C}^{0}_{\frac{\pi}{2} - \theta}$-dominant and $\widehat{{\mathbf B}}$ is NDD, then every $\lambda \in \sigma({\mathbf A})$ satisfies
$$|{\rm arg}(\lambda)| > \theta.$$
\end{theorem}
{\bf Proof.} Let $\lambda = x+iy \in {\mathbb C}$ satisfies $|{\rm arg}(\lambda)| \leq \theta \leq \frac{\pi}{4}.$ Then it satisfies the conditions
 \begin{equation}\label{angl1} x^2 - y^2 \geq 0\end{equation} and
 \begin{equation}\label{angl2} x = |\lambda|\cos({\rm arg}(\lambda)) \geq \cos\theta |\lambda| > 0. \end{equation}
Let us recall that for a given value $\theta \in (0, \frac{\pi}{4}]$, a real $n \times n$ matrix $\widehat{{\mathbf A}}$ is {\it diagonally ${\mathbb C}^{0}_{\frac{\pi}{2} - \theta}$-dominant} if and only if the following inequalities hold:
 \begin{enumerate}
\item[\rm 1.]$\cos \theta|\widehat{a}_{ii}| > \sum_{j\neq i}|\widehat{a}_{ij}| \qquad i = 1, \ \ldots, \ n.$
\item[\rm 2.] $\widehat{a}_{ii} < 0$, $i = 1, \ \ldots, \ n.$
\end{enumerate}
Applying Lemma \ref{1}, we obtain that $\lambda$ does not belong to $\sigma({\mathbf A})$ if and only if ${\mathbf G}(\lambda)$ is nonsingular. One of the conditions sufficient for the non-singularity of ${\mathbf G}(\lambda)$ is its diagonal dominance.
Now let us prove the diagonal dominance of the matrix ${\mathbf G}(\lambda)$ whenever $|{\rm arg}(\lambda)| \leq \theta.$ We have
$$\sum_{i \neq j}|g_{ij}(\lambda)| = \sum_{i \neq j}|\lambda \widehat{a}_{ij} + \widehat{b}_{ij}| \leq |\lambda|\sum_{i \neq j}|\widehat{a}_{ij}| + \sum_{i \neq j}|\widehat{b}_{ij}| < $$
 $$< -|\lambda|\sin(\frac{\pi}{2} - \theta)\widehat{a}_{ii} - \widehat{b}_{ii} = -|\lambda|\cos \theta \widehat{a}_{ii} - \widehat{b}_{ii}\leq -|\lambda|\cos({\rm arg}(\lambda))\widehat{a}_{ii} - \widehat{b}_{ii} = $$
 $$-x\widehat{a}_{ii} - \widehat{b}_{ii} \leq (x^2 - y^2)-x\widehat{a}_{ii} - \widehat{b}_{ii} = |{\rm Re}(\lambda^2- \lambda \widehat{a}_{ii} - \widehat{b}_{ii})| \leq $$
  $$|(\lambda^2- \lambda \widehat{a}_{ii} - \widehat{b}_{ii})| = |g_{ii}(\lambda)|.$$
We have proved that ${\mathbf G}(\lambda)$ is strictly diagonally dominant hence nonsingular for every $\lambda$ with $|{\rm arg}(\lambda)| \leq \theta.$ This fact immediately implies the statement of the theorem. $\square$
\begin{corollary}\label{cor}Let $n = 2k$, $\theta \in (0, \frac{\pi}{4}]$ and an $n \times n$ matrix ${\mathbf A}$ be defined by \eqref{ym}, with ${\mathbf A}_{12} = {\mathbf I}$. If ${\mathbf A}_{11}$ and ${\mathbf A}_{22}$ are diagonally ${\mathbb C}^{0}_{\frac{\pi}{2} - \theta}$-dominant and ${\mathbf A}_{21} - {\mathbf A}_{22}{\mathbf A}_{11}$ is NDD, then every $\lambda \in \sigma({\mathbf A})$ satisfies
$$|{\rm arg}(\lambda)| > \theta.$$
\end{corollary}
{\bf Proof.} For the proof, it is enough to notice, that if ${\mathbf A}_{12} = {\mathbf I}$, Formulas \eqref{A} and \eqref{B} give
$$\widehat{{\mathbf A}} = {\mathbf A}_{11} + {\mathbf A}_{22}; \qquad \widehat{{\mathbf B}} = {\mathbf A}_{21} - {\mathbf A}_{22}{\mathbf A}_{11}.$$
Obviously, if both ${\mathbf A}_{11}$ and ${\mathbf A}_{22}$ are diagonally ${\mathbb C}^{0}_{\frac{\pi}{2} - \theta}$-dominant, then their sum $\widehat{{\mathbf A}}$ is also.
 $\square$

Now we obtain the following sufficient conditions for the asymptotic stability of System \eqref{sysf} and the family of perturbed systems \eqref{sysfper}.
\begin{theorem} Given a linear system of Form \eqref{sysf} with $0 <\gamma \leq \frac{1}{2}$ and $n = 2k$. Let the system matrix $\mathbf A$, represented in Form \eqref{A}, satisfy the following conditions: ${\mathbf A}_{12} = {\mathbf I}$, ${\mathbf A}_{11}$ and ${\mathbf A}_{22}$ are diagonally ${\mathbb C}^{0}_{\frac{\pi}{2}(1 - \gamma)}$-dominant and ${\mathbf A}_{21} - {\mathbf A}_{22}{\mathbf A}_{11}$ is NDD. Then System \eqref{sysf} is asymptotically stable, moreover, each system of the perturbed family \eqref{sysfper}, where
\begin{equation}\label{D} {\mathbf D} = {\rm diag}\{{\mathbf I}, \widehat{{\mathbf D}}\},\end{equation}
 ${\mathbf I}$ is a $k \times k$ identity matrix, $\widehat{{\mathbf D}}$ is a $k \times k$ positive diagonal matrix, is also asymptotically stable.
\end{theorem}
{\bf Proof.} Applying Corollary \ref{cor}, we obtain, that each eigenvalue $\lambda$ of the system matrix $\mathbf A$ satisfies the condition
$$|{\rm arg}(\lambda)| > \theta.$$
This fact immediately implies the asymptotic stability of System \eqref{sysf}.

Now let us consider the matrix ${\mathbf A}_{\mathbf D}$ of the perturbed system \eqref{sysfper}:
$${\mathbf A}_{\mathbf D} = {\mathbf D}{\mathbf A} = \begin{pmatrix}{\mathbf I} & {\mathbf O} \\ {\mathbf O} & \widehat{{\mathbf D}}  \end{pmatrix}\begin{pmatrix} {\mathbf A}_{11} & {\mathbf I} \\
{\mathbf A}_{21} & {\mathbf A}_{22} \end{pmatrix} = \begin{pmatrix} {\mathbf A}_{11} & {\mathbf I} \\
\widehat{{\mathbf D}}{\mathbf A}_{21} & \widehat{{\mathbf D}}{\mathbf A}_{22} \end{pmatrix}. $$
Let us show, that ${\mathbf A}_{\mathbf D}$ satisfies the conditions of Corollary \ref{cor} for every matrix $\mathbf D$ of Form \eqref{D}. Indeed, by Lemma \ref{obv}, if ${\mathbf A}_{22}$ is diagonally ${\mathbb C}^{0}_{\frac{\pi}{2}(1 - \gamma)}$-dominant, then so is $\widehat{{\mathbf D}}{\mathbf A}_{22}$ for every positive diagonal matrix $\widehat{{\mathbf D}}$. Then,
$$\widehat{{\mathbf D}}{\mathbf A}_{21} - (\widehat{{\mathbf D}}{\mathbf A}_{22}){\mathbf A}_{11} = \widehat{{\mathbf D}}({\mathbf A}_{21} - {\mathbf A}_{22}{\mathbf A}_{11}) $$
is NDD whenever ${\mathbf A}_{21} - {\mathbf A}_{22}{\mathbf A}_{11}$ is NDD and $\widehat{{\mathbf D}}$ is positive diagonal. Applying Corollary \ref{cor}, we obtain that every $\lambda \in \sigma({\mathbf A}_{\mathbf D})$ satisfies
$$|{\rm arg}(\lambda)| > \theta.$$ Thus every system of the perturbed family \eqref{sysfper}, is asymptotically stable. $\square$

\subsection{Case $1 \leq\gamma < 2$.} For $1 \leq\gamma < 2$, Condition \eqref{confrac} corresponds to ${\mathfrak D}$-stability of the system matrix $\mathbf A$ with respect to the stability region ${\mathfrak D} = {\mathbb C}^-_{\theta}$, where $\theta = \pi(1 - \frac{\gamma}{2})$.

Let the system matrix $\mathbf A$ be $({\mathbb C}^0_{\theta}, {\mathcal D})$-stable. Then each system of the perturbed family
\begin{equation}\label{sysfper}d^{\gamma}x(t)  = ({\mathbf D}{\mathbf A})x(t), \end{equation}
with $\gamma = \frac{2(\pi - \theta)}{\pi}$, is asymptotically stable for every positive diagonal matrix $\mathbf D$. By Theorem \ref{mainth}, diagonal ${\mathbb C}^0_{\theta}$-dominance of a matrix $\mathbf A$ implies ${\mathbb C}^0_{\theta}$-stability and $({\mathbb C}^0_{\theta}, D)$-stability of $\mathbf A$. Thus we immediately obtain the following result.

\begin{theorem} Given a linear system of the form \eqref{sysf} with $1 \leq\gamma < 2$. Let the system matrix $\mathbf A$ be diagonally ${\mathbb C}^0_{\theta}$-dominant for $\theta = \pi(1 - \frac{\gamma}{2})$. Then System \eqref{sysf} is asymptotically stable, moreover, each system of perturbed family \eqref{sysfper}, where $\mathbf D$ is an arbitrary positive diagonal matrix, is also asymptotically stable.
\end{theorem}

{\bf Proof.} Let $\mathbf A$ be diagonally ${\mathbb C}^0_{\theta}$-dominant. Applying Lemma \ref{obv}, we obtain that ${\mathbf D}{\mathbf A}$ is also diagonally ${\mathbb C}^0_{\theta}$-dominant for any positive diagonal matrix $\mathbf D$. Then, applying Theorem \ref{mmain} to ${\mathbf D}{\mathbf A}$, we complete the proof. $\square$

\section{Numerical examples and simulations}
\bigskip

In this Section we report some numerical examples to illustrate the features
of the presented theoretical approach. All the examples are implemented in
MATLAB 2019a on a genuine Intel Core i7 PC.

\bigskip

\textbf{Example 1 }

The first examples refers to Theorem \ref{diagcrit}.
Let%
\textbf{\ }$\mathbf{A}$\textbf{\ }and\textbf{\ }$\mathbf{B}$\textbf{\ }the
considered real $n\times n$ matrices and $\widetilde{\mathbf{A}}$ be defined
by (8). Let $\mathbf B$ be diagonalizable by an $n \times n$ invertible matrix, say, $\mathbf S$. If $\mathbf{A}$\textbf{\ }and\textbf{\ }$\mathbf{B}$ commute, the
matrices $\mathbf{S}^{-1}\mathbf{AS}$ and $\mathbf{S}^{-1}\mathbf{BS}$ are both
diagonal and eigenvalues of $\mathbf{A}$\textbf{\ }and\textbf{\ }$\mathbf{B,}
$ respectively, appear on the diagonals.

We consider the following commutative matrices

$${\mathbf A} =
\begin{pmatrix}
-0.7143 & -0.3333 \\
-0.1429 & -1.3333%
\end{pmatrix} \qquad {\mathbf B} =
\begin{pmatrix}
-1.4737 & 0.3684 \\
0.1579 & -0.7895%
\end{pmatrix}$$
We can check that the matrices commute by computing $\mathbf{AB=BA}$\textbf{%
\ }$=\mathbf{I.}$

Eigenvalues of $\mathbf{B}$ are $\nu _{1}=-0.713\,03,\nu
_{2}=-1.\,\allowbreak 550\,2\allowbreak $;\ eigenvalues of $\mathbf{A}$ are $%
\mu _{1}=$ $-0.645\,09,$ $\mu _{2}=-1.\,\allowbreak 402\,5\allowbreak $; so
we see that $\mathbf{S}^{-1}\mathbf{AS}$ is diagonally stable. Then $%
\widetilde{\mathbf{A}}$ becomes

$$\widetilde{\mathbf{A}}=\begin{pmatrix}
-0.7143 & -0.3333 & -1.4737 & 0.3684 \\
-0.1429 & -1.3333 & 0.1579 & -0.7895 \\
1 & 0 & 0 & 0 \\
0 & 1 & 0 & 0%
\end{pmatrix}$$
with eigenvalues $\lambda _{1}=$ $-0.322\,55+1.\,\allowbreak
202\,6\allowbreak i,$ $\lambda _{2}=-0.322\,55-1.\,\allowbreak
202\,6\allowbreak i,$ $\lambda _{3}=-0.701\,25+0.470\,39\allowbreak i,$ $%
\lambda _{4}=-0.701\,25-0.470\,39\allowbreak i.=\allowbreak
-0.701\,25-0.470\,39.$ Therefore $\widetilde{\mathbf{A}}$ is stable, as
expected.

\bigskip

\textbf{Example 2}

The second example refers to Lemma \ref{1} and Theorem \ref{14}.  We deal with $n\times
n$ matrix  $\mathbf{A}$ given by \eqref{ym}. Lemma \ref{1} relates to the eigenvalues
of $\mathbf{A}$. \  In \cite{NIS2} the considered block matrix $\mathbf{A}$ is
given by the following blocks%
$$
\mathbf{A}_{11} =
\begin{pmatrix}
-2 & 1.8 \\
1.8 & -2%
\end{pmatrix};
\qquad {\mathbf A}_{12}=
\begin{pmatrix}
-1 & 0.8 \\
-0.8 & -1%
\end{pmatrix};$$
$${\mathbf A}_{21} =
\begin{pmatrix}
1 & 0 \\
0 & 1%
\end{pmatrix}; \qquad {\mathbf A}_{22}=
\begin{pmatrix}
0 & 0 \\
0 & 0%
\end{pmatrix}.$$

In \cite{NIS2}, these matrices allow to show, as a counterexample, that for none of
the theorems presented in that paper  these matrices obey all assumptions,
but still the corresponding matrix $\mathbf{A}$ is stable.

Here,  we easily compute zeros of equation
$$\det \mathbf{G}(\lambda )=0,$$
where $\mathbf{G}(\lambda )=\lambda ^{2}-\widehat{\mathbf{A}}$ $\mathbf{%
\lambda -}\widehat{\mathbf{B}},$ with $\widehat{\mathbf{A}}$ and $\widehat{%
\mathbf{B}}$ are given by \eqref{A} and \eqref{B}, respectively. This means
$$\widehat{\mathbf{A}}=
\begin{pmatrix}
-0.243\,91 & 0.395\,15 \\
0.395\,12 & -3.\,\allowbreak 756\,1%
\end{pmatrix}; \qquad \widehat{\mathbf{B}}=
\begin{pmatrix}
-1.0 & 0.8 \\
-0.8 & -1.0%
\end{pmatrix}
.$$

Since we have $\det \mathbf{G}(\lambda )=\lambda ^{4}+4.0\lambda
^{3}+2.0\lambda^{2}+\allowbreak 4.0\lambda +1.\,\allowbreak 64,$ the
wanted zeros are $\lambda _{1}=$ $-3.\,\allowbreak 499\,7,$ $\lambda
_{2}=-0.470\,88,$ $\lambda _{3}=-0.01471\,9+0.997\,49i,$ $\lambda
_{4}=-0.01471\,9-0.997\,49i.$

As they are the eigenvalues of $\mathbf{A},$ we can conclude that $\mathbf{A}
$ is stable, as expected from \cite{NIS2}.

Now we can use matrices \ $\widehat{\mathbf{A}}$ and $\widehat{\mathbf{B}}$
as a counterexample to show that even if all the hypotheses of Theorem \ref{14}
are not satisfied, still the the thesis happens, since the theorem provides
a sufficient condition only.

Let us assume $\theta =\pi /4=\allowbreak 0.785\,40$; then $%
\widehat{\mathbf{A}}$ is not
 diagonally $\mathbb{C}_{\pi /4}^{0}$-dominant, but $\widehat{\mathbf{B}}$ is NDD.

About $\arg (\lambda )$ with $\lambda \in \sigma (\mathbf{A}),$ we have

\qquad $\arg (-3.\,\allowbreak 499\,7)=\allowbreak 3.\,\allowbreak 141\,6;$

\qquad $\arg (-0.470\,88)=\allowbreak 3.\,\allowbreak 141\,6;$

\qquad $\arg (-0.01471\,9+0.997\,49i)=\allowbreak 1.\,\allowbreak 585\,6;$

\qquad $\arg (-0.01471\,9-0.997\,49i)=\allowbreak -1.\,\allowbreak 585\,6.$

Clearly, for each $\lambda \in \sigma (\mathbf{A})$, $\left\vert \arg
(\lambda )\right\vert >\theta $ \ holds, unexpectedly.

\textbf{Example 3.}

Relating to Corollary \ref{cor}, we introduce the following matrix

$$
{\mathbf A}=
\begin{pmatrix}
{\mathbf A}_{11} & {\mathbf A}_{12} \\
{\mathbf A}_{21} & {\mathbf A}_{22}%
\end{pmatrix},$$
where $${\mathbf A}_{11}=
\begin{pmatrix}
-10 & 3 \\
2 & -10%
\end{pmatrix}; \qquad{\mathbf A}_{12}=
\begin{pmatrix}
1 & 0 \\
0 & 1%
\end{pmatrix};$$ $${\mathbf A}_{21}=
\begin{pmatrix}
20 & 1 \\
2 & 37%
\end{pmatrix}; \qquad {\mathbf A}_{22}=
\begin{pmatrix}
-8 & 4 \\
3 & -8%
\end{pmatrix}.$$

We assume $\theta =\pi /4;\ $\ consequently, ${\mathbf A}_{11}$ and ${\mathbf A}_{22}$ are
 diagonally $\mathbb{C}_{\pi /4}^{0}$-dominant (see Definition
7') since they are NDD\ and the following inequalities hold:\medskip

\qquad $\sin (\pi /4)\ast \left\vert -10\right\vert =\allowbreak
7.\,\allowbreak 07>3$

$\qquad \sin (\pi /4)\ast \left\vert -10\right\vert =\allowbreak
7.\,\allowbreak 07>2$

$\qquad \sin (\pi /4)\ast \left\vert -8\right\vert =\allowbreak
7.\,\allowbreak 07>4$

$\qquad \sin (\pi /4)\ast \left\vert -8\right\vert =\allowbreak
7.\,\allowbreak 07>3\medskip $

Then we have the matrix

$${\mathbf A}_{21}-{\mathbf A}_{22}{\mathbf A}_{11}=
\begin{pmatrix}
-68.0 & 65.0 \\
48.0 & -52.0%
\end{pmatrix},$$
which is NDD.

Therefore all the hypotheses of Corollary \ref{cor} are satisfied. Then we compute
eigenvalues of $\mathbf A$. They all are negative real: \medskip

\qquad $\lambda _{1}=$ $-17.23,$ \ $\lambda _{2}=$ $-0.3017,$ \ $\lambda
_{3}=-11.52,$ \ \ $\lambda _{4}=$ $-6.945$.\medskip

Since the absolute value of the argument of any real number is $\pi $ and $%
\pi >\pi /4,$ the thesis is verified.

\section*{Acknowledgements} The research was supported by the National Science Foundation of China grant number
12050410229.

\end{document}